\documentclass[11pt]{amsart}
\usepackage[T1]{fontenc}
\usepackage[utf8]{inputenc}
\usepackage{lmodern}
\usepackage{amsmath,amssymb,amsthm,mathtools}
\usepackage[margin=1in,headheight=14pt,footskip=30pt]{geometry}
\usepackage{needspace}
\usepackage{microtype}
\usepackage[hidelinks]{hyperref}
\theoremstyle{plain}
\newtheorem{theorem}{Theorem}[section]
\newtheorem{proposition}[theorem]{Proposition}
\newtheorem{lemma}[theorem]{Lemma}
\newtheorem{corollary}[theorem]{Corollary}
\theoremstyle{definition}
\newtheorem{example}[theorem]{Example}
\newtheorem{question}[theorem]{Question}
\theoremstyle{remark}
\newtheorem{remark}[theorem]{Remark}
\numberwithin{equation}{section}
\title[Automorphisms of Chebyshev covering curves]{Automorphisms of Kummer and Artin--Schreier Curves Defined by Chebyshev Polynomials}
\author{Ahmad Kazemifard}
\address{Department of Mathematics, Faculty of Mathematical Sciences and Computer,
Shahid Chamran University of Ahvaz, Ahvaz, Iran}
\email{a.kazemifard@scu.ac.ir}
\author{Saeed Tafazolian}
\address{Departamento de Matem\'atica, Instituto de Matem\'atica, Estat\'istica
e Computa\c{c}\~ao Cient\'ifica (IMECC), Universidade Estadual de Campinas
(UNICAMP), Rua S\'ergio Buarque de Holanda, 651, Cidade Universit\'aria,
Zeferino Vaz, Campinas, SP 13083-859, Brazil}
\email{saeed@unicamp.br}
\hypersetup{pdftitle={Automorphisms of Kummer and Artin--Schreier Curves Defined by Chebyshev Polynomials},pdfauthor={Ahmad Kazemifard and Saeed Tafazolian}}
\begin{document}
\begin{abstract}
Let $\varphi_m$ be the monic Chebyshev polynomial characterized by
$\varphi_m(t+t^{-1})=t^m+t^{-m}$. We study the geometric automorphism
groups of the smooth projective models of $y^n=\varphi_m(x)$ and
$y^p-y=\varphi_m(x)$. We classify the hyperelliptic family in genus
at least two, in characteristic zero and odd characteristic. For every
$m\geq5$, with \(p\nmid m\) in positive characteristic, we determine the projective stabilizer of the roots of
$\varphi_m$ together with infinity. For tame Kummer covers, we describe the
normalizer of the cyclic deck group and compute it explicitly for
$m\geq6$. It is the full group when $m>2n$. If $m\geq3$ and $n>2m$,
the full group has order $2n$ in characteristic zero and in every
characteristic prime to $2nm$. The proof uses the dihedral Galois
closure of the map to the $y$-line and Dickson's subgroup
classification. In odd characteristic, we also classify all
Artin--Schreier curves in the family with genus at least two. For $m\geq3$ prime to $p\geq5$, the
group has order $2p$. In characteristic three, exactly five degrees
are exceptional. We give the groups and their actions, and treat
the quadratic cases separately.
\end{abstract}
\maketitle

\section{Introduction}

The normalized Chebyshev polynomials are characterized by
\begin{equation}\label{eq:intro-Cheb}
 \varphi_m(t+t^{-1})=t^m+t^{-m}.
\end{equation}
They satisfy $\varphi_0=2$, $\varphi_1=X$, and
$\varphi_{m+1}=X\varphi_m-\varphi_{m-1}$. Their coefficients are
integers, so the same normalization can be used in characteristic zero
and in positive characteristic. We consider the smooth projective
models of
\[
 X_{n,m}:\quad y^n=\varphi_m(x),
 \qquad
 A_m:\quad y^p-y=\varphi_m(x).
\]
The constant field $k$ is algebraically closed. In the Kummer case we
assume that its characteristic is zero or prime to $2nm$; the
Artin--Schreier case is considered in characteristic $p>2$ with
$p\nmid m$. These restrictions on $m$ can be removed by the elementary
reductions in Proposition~\ref{prop:Frobenius-reduction}.

The maximality of Chebyshev curves over finite fields has been studied in
\cite{DTT2025,KTT2018,TT2019,TTcorr,TT2026}. Here we determine their
automorphisms over an algebraic closure. For this, we must keep track
of the branch points and their exponents. We must also decide whether
every automorphism preserves the given cover. This holds for a
hyperelliptic curve of genus at least two, since its degree-two map is
unique. For a general Kummer cover it needs a separate proof.

Write $\mathcal Z_m$ for the zero set of $\varphi_m$ and put
\[
 H_m=\operatorname{Stab}_{\operatorname{PGL}_2(k)}(\mathcal Z_m),
 \qquad
 K_m=\operatorname{Stab}_{\operatorname{PGL}_2(k)}
          (\mathcal Z_m\cup\{\infty\}).
\]
Tafazolian and Top \cite{TTprojective} determined $H_m$ for $m>4$.
We recall their result in \eqref{eq:hyp-companion-input}. We determine
$K_m$ for $m\geq6$ in Proposition~\ref{prop:hyp-odd-stabilizer} and
treat $m=5$ separately. The set $\mathcal Z_m\cup\{\infty\}$ occurs
in odd hyperelliptic degree. It also occurs when infinity has the same
Kummer exponent as the finite roots. Our proof adapts the
logarithmic-differential argument of \cite[Section~4]{TTprojective}.
The branch form now has degree $m+1$, and the exceptional case is
$2m+1=p^e$. The plane case $n=m$ was treated in \cite{TTplane}.
We use these earlier results below.

Theorem~\ref{thm:hyp-classification} gives the full hyperelliptic
classification. It includes the degree-five curve and all exceptions
in odd characteristic. In characteristic zero the answer is
\[
 \operatorname{Aut}(X_{2,m})\cong
 \begin{cases}
 C_2\times C_2,&m\geq6\text{ even},\\
 D_4,&m=5,\\
 C_4,&m\geq7\text{ odd}.
 \end{cases}
\]
Here and throughout, $C_a$ denotes a cyclic group of order $a$, and
$D_a$ a dihedral group of order $2a$. In particular, $D_4$ has order
eight. The lift of $x\mapsto-x$ has order four in odd degree, which
explains the cyclic group in the last row.

For a Kummer curve, the branch exponents are taken modulo $n$.
They force the cyclic deck group to be central in its normalizer.
Theorem~\ref{thm:kummer-normalizer} describes this normalizer, and
Theorem~\ref{thm:explicit-normalizers} gives its group structure for
$m\geq6$. It is the full automorphism group when $m>2n$.
This follows from the strong-branching criterion of
Kontogeorgis~\cite{Kontogeorgis1999}. We include a proof that also
applies when $n$ is composite.

We then consider the map to the $y$-line. Its dihedral Galois closure
allows us to determine the full group for $m\geq3$ and $n>2m$;
see Theorem~\ref{thm:full-Cheb-high-n-sharpened}. When $m$ is odd,
$m\mid n$ and $2n=p^r+1$, the branch set is projectively equivalent
to $\mathbb P^1(\mathbb F_{p^r})$. In this case we use Dickson's
classification together with Lemma~\ref{lem:Kummer-no-inversion},
which excludes lifts of inversions. For every $m\geq3$ and $n>2m$,
the answer is
\[
 \mathcal A_{n,m}\cong
 \begin{cases}
 C_n\times C_2,&m\text{ even},\\
 C_{2n},&m\text{ odd}.
 \end{cases}
\]
These statements hold in characteristic zero and in every
characteristic prime to $2nm$. The case $m=2$ is classified in
Corollary~\ref{cor:quadratic-Cheb-full}. Section~\ref{sec:remaining}
gives examples where the deck group is not normal and states the
remaining question.

For Artin--Schreier curves, we obtain a complete classification in
genus at least two and odd characteristic. If $m\geq3$ and $p\geq5$,
the group is $C_p\times C_2$ for even $m$ and $D_p$ for odd $m$.
In characteristic three, the exceptional degrees are $4,5,7,8,14$;
see Theorem~\ref{thm:AS-full-classification}. We give the groups and
explicit formulas for their actions. We first use Stichtenoth's
theorem, in the form recalled by Lehr and
Matignon~\cite{LehrMatignon2005}, to show that every automorphism fixes
the point above infinity. The rest of the proof is a comparison of
coefficients. Proposition~\ref{prop:AS-degree-two} treats degree two.

\section*{Use of AI tools}
We used ChatGPT (OpenAI)  to assist with
mathematical arguments, symbolic and numerical checks, and the
English text. We take responsibility for the final manuscript and
all mathematical claims.

\section{Removing powers of the characteristic}

All curves in this paper are smooth projective models of their
displayed affine equations. This convention is necessary when a
polynomial has repeated roots. The following reduction also explains
which degree is relevant for the genus and the exceptional cases.

\begin{proposition}\label{prop:Frobenius-reduction}
Suppose $\operatorname{char}(k)=p>0$ and write
$m=p^s m_0$, with $p\nmid m_0$. Then
\[
 \varphi_m(X)=\varphi_{m_0}(X)^{p^s}.
\]
If $p\nmid n$, the smooth projective model of
$y^n=\varphi_m(x)$ is isomorphic to $X_{n,m_0}$ by an isomorphism
preserving $x$ and the cyclic deck subgroup. The smooth projective
model of $y^p-y=\varphi_m(x)$ is likewise isomorphic to $A_{m_0}$.
\end{proposition}
\begin{proof}
The polynomial identity follows from \eqref{eq:intro-Cheb} and
Frobenius. Put $q=p^s$ and $f=\varphi_{m_0}(x)$. Choose integers
$a,b$ with $aq+bn=1$. In the Kummer function field define
\[
 z=y^a f^b.
\]
Then $z^n=f$ and $z^q=y$. Thus $k(x,y)=k(x,z)$, and the deck
generator $y\mapsto\zeta_ny$ becomes
$z\mapsto\zeta_n^a z$, another generator of the same cyclic subgroup.
Negative exponents give rational functions and cause no difficulty on
the smooth projective models.

For the Artin--Schreier equation put
\[
 z=y-\sum_{j=0}^{s-1}f^{p^j},
\]
where the sum is zero when $s=0$. The identity
$(\sum_{j=0}^{s-1}f^{p^j})^p-\sum_{j=0}^{s-1}f^{p^j}=f^q-f$
gives $z^p-z=f$, as required.
\end{proof}

All degree conditions below apply to the part of $m$ prime to the
characteristic. In odd characteristic, the Kummer genus is
\[
 \frac{(n-1)(m_0-1)-\gcd(n,m_0)+1}{2},
\]
and the equal-weight condition at infinity is $n\mid(m_0+1)$.
We henceforth assume $p\nmid m$ in positive characteristic.

We work in characteristic zero or odd characteristic. The
Artin--Schreier reduction also works in characteristic two. Our later
arguments, however, use $x\mapsto-x$ and the distinct critical values
$2,-2$. The Kummer covers are tame, so $p\nmid n$. Throughout the
paper, $\operatorname{Aut}$ denotes the full geometric group.

\section{The hyperelliptic curves}
\label{sec:hyp}

Throughout this section, $k$ is algebraically closed of characteristic zero
or of odd characteristic $p$, and $\varphi_m$ denotes the normalized Chebyshev
polynomial characterized by
\[
  \varphi_m(t+t^{-1})=t^m+t^{-m}.
\]
In positive characteristic we assume $p\nmid m$. Let $X_m$ be the smooth
projective curve with function field $k(x,y)$, where
\[
  y^2=\varphi_m(x).
\]
Thus $X_m=X_{2,m}$ in the notation of the introduction.
We determine its full automorphism group for $m\geq 5$. We write $C_a$ for
the cyclic group of order $a$ and $D_a$ for the dihedral group of order $2a$.
In particular, $D_4$ has order eight.

For an even integer $N\geq 6$, define the group
\begin{equation}
\label{eq:hyp-GN}
 \mathcal G_N=
 \left\langle r,h,t\;\middle|\;
 \begin{array}{c}
 r^N=h^2=t^2=1,\quad [r,h]=[t,h]=1,\\
 trt^{-1}=hr^{-1}
 \end{array}
 \right\rangle.
\end{equation}
Thus $\mathcal G_N=(C_N\times C_2)\rtimes C_2$, where the last factor fixes
$h$ and sends $r$ to $hr^{-1}$. This action has order two because $N$ is
even. Thus $|\mathcal G_N|=4N$. A central extension of a dihedral
group need not be dihedral, so we give the presentation explicitly.

\begin{theorem}
\label{thm:hyp-classification}
Suppose $m\geq 5$, and retain the assumptions on $k$ and $m$ above.
For even $m$, one has
\begin{equation}
\label{eq:hyp-even-classification}
 \operatorname{Aut}(X_m)\cong
 \begin{cases}
  \mathcal G_m,& \text{if }2m=p^e+1\text{ for some }e\geq1,\\
  C_2\times S_3,& \text{if }4m=p^e+1\text{ for some }e\geq1,\\
  C_2\times C_2,& \text{otherwise}.
 \end{cases}
\end{equation}
For odd $m$, one has
\begin{equation}
\label{eq:hyp-odd-classification}
 \operatorname{Aut}(X_m)\cong
 \begin{cases}
  \mathcal G_{m+1},& \text{if }2m+1=p^e\text{ for some }e\geq1,\\
  D_4,& \text{if }m=5\text{ and }\operatorname{char}(k)\ne11,\\
  C_4,& \text{if }m\geq7\text{ and no such }p^e\text{ exists}.
 \end{cases}
\end{equation}
Conditions involving $p^e$ apply only in positive characteristic. In
particular, in characteristic zero the groups are $C_2\times C_2$ for even
$m$, $D_4$ for $m=5$, and $C_4$ for odd $m\geq7$.
\end{theorem}

For even $m$, we lift the known stabilizer of the zeros of $\varphi_m$
and identify the central extension. For odd $m$, we compute the branch
stabilizer below. The degree-five curve needs a separate calculation.

\subsection{The branch divisor and the central extension}
\label{subsec:hyp-branch}

Put $f=\varphi_m$. The polynomial identities
\begin{equation}
\label{eq:hyp-cheb-identities}
 (x^2-4)f'(x)^2=m^2\bigl(f(x)^2-4\bigr),
 \qquad
 (x^2-4)f''(x)+xf'(x)-m^2f(x)=0
\end{equation}
hold over $k$. They follow by differentiating the defining identity for
$\varphi_m$ and are polynomial identities over the integers. Since
$f(2)=2$ and $f(-2)=2(-1)^m$, the points $2$ and $-2$ are not roots.
The first identity shows that every root is simple. Moreover,
\begin{equation}
\label{eq:hyp-leading-coefficients}
 f(x)=x^m-mx^{m-2}+\text{terms of lower degree},
 \qquad f(-x)=(-1)^m f(x).
\end{equation}

Let $\mathcal Z_m\subset\mathbb P^1(k)$ be the set of zeros of $f$.
The branch set of $x:X_m\longrightarrow\mathbb P^1$ is
\begin{equation}
\label{eq:hyp-branch-set}
 B_m=
 \begin{cases}
  \mathcal Z_m,&m\text{ even},\\
  \mathcal Z_m\cup\{\infty\},&m\text{ odd}.
 \end{cases}
\end{equation}
Consequently $g(X_m)=\lfloor(m-1)/2\rfloor$. In particular, $g(X_m)\geq2$
under our degree assumption.

The hyperelliptic involution is $h:(x,y)\mapsto(x,-y)$.
The degree-two map is unique up to an automorphism of $\mathbb P^1$.
Indeed, two distinct degree-two rational subfields would give a
birational map to a curve of bidegree $(2,2)$ in
$\mathbb P^1\times\mathbb P^1$. That curve has arithmetic genus one,
contrary to $g(X_m)\geq2$. Thus every
automorphism of $X_m$ descends to a projective transformation preserving
$B_m$. Conversely, for any such transformation $M$, the rational function
$f(M(x))/f(x)$ has even divisor. It is a square in $k(x)$, since $k$ is
algebraically closed. Choosing either square root gives a lift of $M$.
It follows that
\begin{equation}
\label{eq:hyp-exact-sequence}
 1\longrightarrow\langle h\rangle
 \longrightarrow\operatorname{Aut}(X_m)
 \longrightarrow \overline G_m
 \longrightarrow1,
 \qquad
 \overline G_m=\operatorname{Stab}_{\operatorname{PGL}_2(k)}(B_m),
\end{equation}
is exact, with central kernel of order two. In odd degree, the branch
set includes infinity.

We record the realization of the group in
\eqref{eq:hyp-GN}. Assume $N$ is even, $N\geq6$, and
$\operatorname{char}(k)\nmid2N$. On the smooth curve
\begin{equation}
\label{eq:hyp-binomial-model}
 W_N:\quad w^2=z^N+1,
\end{equation}
choose a primitive $N$-th root of unity $\zeta_N$. The maps
\begin{equation}
\label{eq:hyp-GN-actions}
 r(z,w)=(\zeta_Nz,w),\qquad
 h(z,w)=(z,-w),\qquad
 t(z,w)=\left(z^{-1},\frac{w}{z^{N/2}}\right)
\end{equation}
are automorphisms. They satisfy precisely the relations displayed in
\eqref{eq:hyp-GN}: for example,
$\zeta_N^{N/2}=-1$ gives $trt^{-1}=hr^{-1}$. Every word has a normal form
$r^a h^b t^c$, with $0\leq a<N$ and $b,c\in\{0,1\}$; their actions are
distinct, so this subgroup has order $4N$. Its reduced image is $D_N$.
If the branch stabilizer of $W_N$ is $D_N$, this subgroup is the full
automorphism group.

\subsection{The augmented zero set}
\label{subsec:hyp-odd}

We adapt the calculation in \cite[Section~4]{TTprojective} to the set
of zeros together with infinity. We give the details below.
In this subsection $m\geq6$ may have either parity,
$B=\mathcal Z_m\cup\{\infty\}$, $K=\operatorname{Stab}(B)$, and
$\sigma(x)=-x$.

\begin{lemma}
\label{lem:hyp-infinity}
For every $m\geq2$, the stabilizer of infinity in $K$ is
$K_\infty=\{1,\sigma\}$. In particular, $|K|\leq2(m+1)$.
If $K\ne\{1,\sigma\}$, then $K$ contains an involution distinct from
$\sigma$.
\end{lemma}

\begin{proof}
If $x\mapsto ax+b$ preserves $\mathcal Z_m$, monicity gives
$f(ax+b)=a^m f(x)$. Comparison of the coefficient of $x^{m-1}$, followed
by that of $x^{m-2}$, gives $b=0$ and $a^2=1$, respectively, because
$m\ne0$ in $k$. Thus $K_\infty=\{1,\sigma\}$. Since $K$ acts faithfully
on the finite set $B$, it is finite; the order bound follows from the
orbit--stabilizer formula.

Suppose $\sigma$ is the unique involution of $K$. It is then central.
The centralizer of $\sigma$ in $\operatorname{PGL}_2(k)$ consists of
the maps $x\mapsto ax$ and $x\mapsto a/x$. Every map of the latter type is
an involution different from $\sigma$, so it cannot occur in $K$.
The remaining maps fix infinity and hence belong to $\{1,\sigma\}$.
This proves the last assertion.
\end{proof}

\begin{lemma}
\label{lem:hyp-differential}
If $m\geq6$ and $K\ne\{1,\sigma\}$, then $k$ has positive
characteristic $p$ dividing $2m+1$. Moreover, $K$ contains a nonidentity
element fixing both $2$ and $-2$ individually.
\end{lemma}

\begin{proof}
Choose an involution $\tau\ne\sigma$ in $K$ and put $M=\sigma\tau$.
Then $M\notin\{1,\sigma\}$ and
$\sigma M\sigma=M^{-1}$. By Lemma~\ref{lem:hyp-infinity}, $M$ is
nonaffine. Write $u=M(\infty)$. Thus $u\in\mathcal Z_m$, and the
reversibility relation gives
$M^{-1}(\infty)=-u$. Consequently
\begin{equation}
\label{eq:hyp-M-normal-form}
 M(x)=u+\frac{s}{x+u},\qquad s\ne0.
\end{equation}
The case $u=0$ is allowed. In every case $u\ne\pm2$.

The homogeneous form defining $B$ has degree $m+1$, and preservation of its
simple zero set gives
\begin{equation}
\label{eq:hyp-branch-identity}
 (x+u)^{m+1}f(M(x))=\lambda f(x)
 \quad\text{for some }\lambda\in k^*.
\end{equation}
Set
\[
 A(x)=-\frac{x}{x^2-4},\qquad \omega=A(x)\,dx.
\]
At any root $\beta$ of $f$, the second identity in
\eqref{eq:hyp-cheb-identities} gives
$f''(\beta)/f'(\beta)=A(\beta)$. If $\beta\ne-u$, then $M(\beta)$ is
also a finite root. Differentiate
\eqref{eq:hyp-branch-identity} twice and divide its second derivative by
its first derivative at $\beta$. Since
$M''/M'=-2/(x+u)$, this yields
\begin{equation}
\label{eq:hyp-root-differential}
 A(\beta)=\frac{2m}{\beta+u}+M'(\beta)A(M(\beta)).
\end{equation}

We now extend this equality to an identity of rational functions. Put $q=x+u$ and
$R_0=(uq+s)^2-4q^2$. With denominator $(x^2-4)qR_0$, the numerator of
$A(x)-2m/q-M'(x)A(M(x))$ is
\begin{equation}
\label{eq:hyp-numerator}
 -xqR_0-2m(x^2-4)R_0-s(x^2-4)(uq+s),
\end{equation}
which has degree at most four. Its denominator is nonzero at each of the
$m-1$ roots $\beta\ne-u$: neither $\beta$ nor $M(\beta)$ is $\pm2$.
Since $m-1>4$, \eqref{eq:hyp-root-differential} forces
\eqref{eq:hyp-numerator} to vanish identically. Thus
\begin{equation}
\label{eq:hyp-one-form}
 \omega=\frac{2m\,dx}{x+u}+M^*\omega.
\end{equation}

The residues of $\omega$ at $2,-2,\infty$ are
$-1/2,-1/2,1$, respectively. Its left-hand side in
\eqref{eq:hyp-one-form} is regular at $-u$. On the right-hand side,
$M^*\omega$ has residue $1$ there, since $M(-u)=\infty$, and the first
term has residue $2m$. Therefore $2m+1=0$ in $k$. This is impossible in
characteristic zero and gives the asserted divisibility in positive
characteristic.

After the poles at $-u$ cancel, the remaining finite poles of the
right-hand side are exactly $M^{-1}(\{2,-2\})$. These points are finite
because $M(\infty)=u\ne\pm2$, and their residues are nonzero. Equality
with the left-hand side forces
$M^{-1}(\{2,-2\})=\{2,-2\}$. Hence $M$ permutes $2$ and $-2$.
Either $M$ or $\sigma M$ fixes both points individually; that element is
nonidentity because $M\notin\{1,\sigma\}$.
\end{proof}

We now make the resulting scaling condition explicit. Define
\begin{equation}
\label{eq:hyp-PQ}
 P_m(z)=\sum_{j=0}^{m}\binom{2m}{2j}z^j,
 \qquad Q_m(z)=(1-z)P_m(z).
\end{equation}
The coordinate $z=(x-2)/(x+2)$ gives
\begin{equation}
\label{eq:hyp-coordinate}
 (1-z)^m f\!\left(2\frac{1+z}{1-z}\right)=2P_m(z),
 \qquad
 (1-z)^{m+1} f\!\left(2\frac{1+z}{1-z}\right)=2Q_m(z).
\end{equation}
One can verify the first identity by substituting
$t=(1+T)/(1-T)$ and $z=T^2$ into the defining identity for $\varphi_m$.
The polynomial $Q_m$ has degree $m+1$, constant coefficient one, and
simple roots equal to the image of $B$. In particular, $z=1$ is the
image of the added point at infinity. An element fixing $2$ and $-2$
individually becomes a scaling $z\mapsto az$. If it preserves $B$, then
\begin{equation}
\label{eq:hyp-scaling}
 Q_m(az)=Q_m(z),
\end{equation}
where the multiplier on the right is one because $Q_m(0)=1$.

\begin{lemma}
\label{lem:hyp-Lucas}
Assume $m\geq2$, $p\mid2m+1$, and $p$ is odd. The polynomial $Q_m$
admits a scaling as in \eqref{eq:hyp-scaling} with $a\ne1$ if and only if
$2m+1$ is a power of $p$. In that case
\begin{equation}
\label{eq:hyp-Q-Frobenius}
 Q_m(z)=1-z^{m+1}.
\end{equation}
\end{lemma}

\begin{proof}
Write $2m+1=p^eL$ with $e\geq1$ and $p\nmid L$. Both $L$ and $p^e$
are odd. From \eqref{eq:hyp-PQ} one obtains the exact identity
\[
 2Q_m(T^2)=(1-T)(1+T)^{2m+1}+(1+T)(1-T)^{2m+1}.
\]
Frobenius and the binomial theorem give
\begin{equation}
\label{eq:hyp-Q-Lucas-expansion}
 Q_m(T^2)=
 \sum_{\substack{0\leq\ell\leq L\\\ell\ \mathrm{even}}}
 \binom{L}{\ell}T^{p^e\ell}
 -
 \sum_{\substack{0\leq\ell\leq L\\\ell\ \mathrm{odd}}}
 \binom{L}{\ell}T^{p^e\ell+1}.
\end{equation}
The exponents from the two sums cannot coincide, since such a coincidence
would imply that $p^e$ divides $1$.

Write $L=\sum_{i\geq0}L_ip^i$ in base $p$. We use the elementary
binomial-coefficient rule
\[
 \binom{L}{\ell}\equiv\prod_i\binom{L_i}{\ell_i}\pmod p,
 \qquad \ell=\sum_i\ell_ip^i.
\]
It follows directly by comparing coefficients in
$(1+T)^L=\prod_i(1+T^{p^i})^{L_i}$ over $k$.
If some $L_i\geq2$, choose $\ell=2p^i$. The corresponding nonzero
coefficient in \eqref{eq:hyp-Q-Lucas-expansion} occurs at
$z^{p^{e+i}}$. Hence \eqref{eq:hyp-scaling} implies
$a^{p^{e+i}}=1$, forcing $a=1$ in characteristic $p$. Thus all the
digits $L_i$ are zero or one. As $p\nmid L$, we have $L_0=1$.

If some higher digit $L_i$ is one, the three choices
$\ell=1,p^i,1+p^i$ give nonzero coefficients at the respective powers
\[
 E_1=\frac{p^e+1}{2},\qquad
 E_2=\frac{p^{e+i}+1}{2},\qquad
 E_3=\frac{p^e+p^{e+i}}{2}
\]
of $z$. Since $E_1+E_2-E_3=1$, the equations
$a^{E_1}=a^{E_2}=a^{E_3}=1$ again imply $a=1$.
Consequently a nontrivial scaling forces $L=1$.
Conversely, for $L=1$, \eqref{eq:hyp-Q-Lucas-expansion} reads
$Q_m(T^2)=1-T^{p^e+1}$, proving
\eqref{eq:hyp-Q-Frobenius}. Its zeros admit all scalings by
$(m+1)$-st roots of unity; these include nonidentity scalings because
$p\nmid m+1$.
\end{proof}

\begin{proposition}
\label{prop:hyp-odd-stabilizer}
For every $m\geq6$, let $K_m=\operatorname{Stab}_{\operatorname{PGL}_2(k)}(\mathcal Z_m\cup\{\infty\})$. Then
\[
 K_m\cong
 \begin{cases}
  D_{m+1},&2m+1=p^e\text{ for some }e\geq1,\\
  C_2,&\text{otherwise}.
 \end{cases}
\]
In the second case $K_m=\{x\mapsto x,x\mapsto-x\}$.
\end{proposition}

\begin{proof}
If $K_m$ is larger than $\{1,\sigma\}$,
Lemma~\ref{lem:hyp-differential} gives $p\mid2m+1$ and a nontrivial
scaling in the $z$ coordinate. Lemma~\ref{lem:hyp-Lucas} therefore
gives $2m+1=p^e$.

Conversely, in that Frobenius case, \eqref{eq:hyp-Q-Frobenius} identifies
the augmented set $B$ with $\mu_{m+1}$. Its stabilizer contains the rotations
$z\mapsto\zeta z$ for $\zeta^{m+1}=1$ and the inversion
$z\mapsto z^{-1}$, a dihedral subgroup of order $2(m+1)$.
Lemma~\ref{lem:hyp-infinity} already gives $|K_m|\leq2(m+1)$, so this
subgroup is the full stabilizer.
\end{proof}

We now return to odd $m\geq7$, for which $\overline G_m=K_m$.
The automorphism groups follow from
\eqref{eq:hyp-exact-sequence}. In the nonexceptional case, the map
\begin{equation}
\label{eq:hyp-order-four}
 \rho(x,y)=(-x,iy),\qquad i^2=-1,
\end{equation}
has order four and square $h$. Since $|K_m|=2$, it generates the full
group. In the exceptional case, setting
\[
 x=2\frac{1+z}{1-z},\qquad
 w=\frac{(1-z)^{(m+1)/2}}{\sqrt2}\,y
\]
gives $w^2=1-z^{m+1}$. Rescaling $z$ identifies this curve with
$W_{m+1}$ from \eqref{eq:hyp-binomial-model}. Its reduced stabilizer is
$D_{m+1}$, so \eqref{eq:hyp-GN-actions} and
\eqref{eq:hyp-exact-sequence} give the full group $\mathcal G_{m+1}$.

\subsection{The genus-two curve of degree five}
\label{subsec:hyp-five}

When $m=5$, the numerator in Lemma~\ref{lem:hyp-differential} has only
four known zeros. Its degree is at most four, so that argument does
not apply. In this degree the group is larger even in characteristic zero.

\begin{proposition}
\label{prop:hyp-five}
Assume $\operatorname{char}(k)\ne2,5$. For
$X_5:y^2=x^5-5x^3+5x$, the reduced group and full automorphism group are
\[
 \bigl(K_5,\operatorname{Aut}(X_5)\bigr)\cong
 \begin{cases}
  (D_6,\mathcal G_6),&\operatorname{char}(k)=11,\\
  (C_2\times C_2,D_4),&\operatorname{char}(k)\ne11.
 \end{cases}
\]
These are all cases, including characteristic three.
\end{proposition}

\begin{proof}
Write $f=x^5-5x^3+5x$ and choose $a\in k$ with $a^2=5$. The maps
\[
 \sigma(x)=-x,\qquad \iota(x)=\frac{a}{x}
\]
preserve $B=\mathcal Z_5\cup\{\infty\}$ and generate a Klein four
group $V$. Indeed,
$f(a/x)=5a f(x)/x^6$. The involution $\sigma$ fixes precisely the two
branch points $0,\infty$. Neither $\iota$ nor $\sigma\iota$ fixes a
branch point. To see the latter assertion, a finite nonzero root satisfies
$x^4-5x^2+5=0$, whereas a fixed point of either involution satisfies
$x^2=a$ or $x^2=-a$. These conditions would imply $a=2$ or $a=-2$,
contradicting $a^2=5$. Thus $V$ has two orbits on $B$, of lengths two
and four.

By Lemma~\ref{lem:hyp-infinity}, $K_\infty=\{1,\sigma\}$. If
$K=\operatorname{Stab}(B)$ is larger than $V$, its orbit of infinity
contains all six points. Indeed, it contains $\{0,\infty\}$, but an
orbit of length two would give $|K|=4$. It therefore meets the other
$V$-orbit and includes that entire orbit. Hence in this
case $|K|=12$; otherwise $K=V$ and $|K|=4$.

Suppose $K>V$. If its only involutions were the three in $V$, then
$\sigma$, being the only one with two fixed branch points, would be
fixed by conjugation and hence central in $K$. The centralizer of
$\sigma$ consists of $ax$ and $a/x$, which preserve
$\{0,\infty\}$. This contradicts transitivity on $B$. There is
therefore an involution $\tau\notin V$. Put $M=\sigma\tau$.
Then $M\notin V$ and $\sigma M\sigma=M^{-1}$. Its value
$u=M(\infty)$ is a finite nonzero root of $f$: a map sending infinity
to infinity belongs to $K_\infty\subset V$, and a map sending it to
zero belongs to $\iota K_\infty\subset V$. As before,
\[
 M(x)=u+\frac{s}{x+u},\qquad s\ne0.
\]
Set $v=u^2$ and $c=f'(u)$. Then
\begin{equation}
\label{eq:hyp-five-vc}
 v^2-5v+5=0,\qquad c=10(v-2)\ne0.
\end{equation}

We compare the branch forms directly. With $q=x+u$, the identity
$q^6 f(u+s/q)=\lambda f(q-u)$ has leading coefficient
$\lambda=sc$. Dividing by $sq$ gives
\begin{align}
\label{eq:hyp-five-comparison}
 &c q^4+5u(2v-3)s q^3+5(2v-1)s^2 q^2+5u s^3q+s^4\notag\\
 &\hspace{15mm}=
 c\bigl(q^4-5u q^3+5(2v-1)q^2-5u(2v-3)q+c\bigr).
\end{align}
This is a polynomial expansion, so it is also valid in characteristic
three. Since $5u\ne0$, comparison of the coefficients of $q^3,q^2,q$
gives, respectively,
\begin{equation}
\label{eq:hyp-five-three-equations}
 (2v-3)s=-c,\qquad
 (2v-1)(s^2-c)=0,\qquad
 s^3=-(2v-3)c.
\end{equation}
The first and third equations, using $s\ne0$, imply
$s^2=(2v-3)^2$. By \eqref{eq:hyp-five-vc},
\[
 (2v-3)^2-c=9-2v.
\]
The middle equation in \eqref{eq:hyp-five-three-equations} consequently
becomes
\[
 (2v-1)(9-2v)=0.
\]
Thus $v=1/2$ or $v=9/2$. Substitution of either value into
$v^2-5v+5=0$ gives $11/4=0$. It follows that $K>V$ can occur only in
characteristic eleven.

Conversely, in characteristic eleven the identity
$2m+1=11$ and \eqref{eq:hyp-Q-Frobenius} give
$Q_5(z)=1-z^6$. Thus $B$ is projectively equivalent to $\mu_6$, whose
dihedral subgroup has order twelve. The established bound $|K|\leq12$
proves $K=D_6$, and the same change of variables as above identifies
$X_5$ with $W_6$. Hence its full group is $\mathcal G_6$.

In every other allowed characteristic, $K=V$ and
$|\operatorname{Aut}(X_5)|=8$. To determine the extension, choose
$b\in k$ with $b^4=5$ and rescale $x=bU$ and $y=c_0Y$, where
$c_0^2=b^5$. The equation becomes
\[
 Y^2=U^5-b^2U^3+U.
\]
The maps
\[
 \rho(U,Y)=(-U,iY),\qquad
 t(U,Y)=\left(U^{-1},\frac{Y}{U^3}\right),\qquad i^2=-1,
\]
satisfy $\rho^4=t^2=1$ and $t\rho t^{-1}=\rho^{-1}$.
Their reduced images generate $V$, and $\rho^2=h$, so they generate
a subgroup of order eight. It is therefore the full group, isomorphic
to $D_4$.
\end{proof}

\subsection{Even degree}
\label{subsec:hyp-even}

For even $m\geq6$, the branch set is $\mathcal Z_m$.
We use the following result from \cite{TTprojective}. Under our
assumptions on $k$, it holds for every $m>4$ prime to the characteristic:
\begin{equation}
\label{eq:hyp-companion-input}
 H_m=\operatorname{Stab}_{\operatorname{PGL}_2(k)}(\mathcal Z_m)\cong
 \begin{cases}
  D_m,&2m=p^e+1,\\
  S_3,&4m=p^e+1,\\
  C_2,&\text{otherwise},
 \end{cases}
\end{equation}
where the generic $C_2$ is generated by $x\mapsto-x$; in characteristic
zero only the last case occurs
\cite[Theorem~1.1, Remark~5.1 and Propositions~5.2--5.3]{TTprojective}.
The two exceptional equalities cannot hold simultaneously: they would
give $2p^a+1=p^b$ for some positive $a,b$, which is impossible modulo
$p$.

In the generic case, the lift $(x,y)\mapsto(-x,y)$ has order two and
commutes with $h$. The exact sequence
\eqref{eq:hyp-exact-sequence} therefore gives $C_2\times C_2$.

If $2m=p^e+1$, Frobenius gives
$P_m(z)=1+z^m$: in $(1+T)^{p^e+1}$, only the even exponents $0$ and
$p^e+1$ occur. The change of variables
\[
 x=2\frac{1+z}{1-z},\qquad
 w=\frac{(1-z)^{m/2}}{\sqrt2}\,y
\]
identifies $X_m$ with $W_m$. Its reduced group is $D_m$ by
\eqref{eq:hyp-companion-input}, and hence its full group is
$\mathcal G_m$.

Finally, suppose $4m=p^e+1$, so that the reduced group is $S_3$.
The lift $s:(x,y)\mapsto(-x,y)$ is an involution above a transposition.
Choose a lift of a reduced element of order three. It has order three
or six, and in the latter case multiplication by $h$ gives a lift
$r$ of order three. The conjugate $srs^{-1}$ is a lift of the inverse
reduced element and has order three. Of its two lifts, $r^{-1}$ has
order three and $hr^{-1}$ has order six; therefore
$srs^{-1}=r^{-1}$. Thus $r,s$ generate a complement isomorphic to
$S_3$ to the central subgroup $\langle h\rangle$. This proves
$\operatorname{Aut}(X_m)\cong C_2\times S_3$ and completes the proof
of Theorem~\ref{thm:hyp-classification}.

\begin{remark}
\label{rem:hyp-scope}
The assumption $m\geq5$ is equivalent to genus at least two. For $m=1,2$
the smooth curve has genus zero, and its geometric automorphism group
is $\operatorname{PGL}_2(k)$. For $m=3,4$ it has genus one, and the full
group contains the group of translations of the elliptic curve after
choosing an origin. These full groups are infinite. For a genus-one curve, the subgroup
fixing a chosen point is finite.
\end{remark}

\section{The deck normalizer of a Chebyshev--Kummer curve}

Let $k$ be algebraically closed of characteristic zero or of characteristic
$p>0$ with $p\nmid 2nm$. Let $n,m\geq2$, and let $X_{n,m}$ denote the smooth
projective curve with function field
\[
 F=k(x,y),\qquad y^n=\varphi_m(x),
 \qquad \varphi_m(t+t^{-1})=t^m+t^{-m}.
\]
Put $\delta=\gcd(n,m)$, let $\mathcal Z_m$ be the zero set of $\varphi_m$, and set
\[
 N=\langle\nu\rangle\cong C_n,\qquad
 \nu(x,y)=(x,\zeta_n y),\qquad
 M=\operatorname{N}_{\operatorname{Aut}(X_{n,m})}(N).
\]
We assume $g(X_{n,m})\geq2$ whenever the full automorphism group is used.
We then study when $N$ is normal, or equivalently when
$M=\operatorname{Aut}(X_{n,m})$.

\begin{proposition}\label{prop:kummer-genus}
The polynomial $\varphi_m$ is separable. The extension $F/k(x)$ is cyclic of
degree $n$, each of the $m$ finite zeros of $\varphi_m$ is totally ramified,
and there are $\delta$ points above infinity, each of ramification index
$n/\delta$. Consequently
\[
 g(X_{n,m})
 =\frac{(n-1)(m-1)-\delta+1}{2}.
\]
\end{proposition}
\begin{proof}
The identity defining $\varphi_m$ shows that its zeros are the images
$t+t^{-1}$ of the roots of $t^{2m}+1$, with $t$ and $t^{-1}$ identified.
The roots of $t^{2m}+1$ are simple and none is $1$ or $-1$. Thus there are
exactly $m$ distinct images. At any finite zero of $\varphi_m$ its valuation
is one, so $Y^n-\varphi_m(x)$ is Eisenstein in the corresponding discrete
valuation ring. It is therefore irreducible, and the extension is cyclic
because $k$ contains the $n$th roots of unity. The tame Kummer ramification
formula \cite[Proposition~3.7.3]{Stichtenoth} applied to the valuations
$1$ and $-m$ gives the asserted indices.
Riemann--Hurwitz now gives
\[
 2g-2=-2n+m(n-1)+(n-\delta)
      =nm-n-m-\delta.
\]
\end{proof}

For $P\in\mathbb P^1(k)$, write
\[
 w(P)=\operatorname{ord}_P(\varphi_m)\pmod n.
\]
Thus $w(P)=1$ for $P\in \mathcal Z_m$, $w(\infty)=-m$, and $w(P)=0$ elsewhere.
These weights lie in $\mathbb Z/n\mathbb Z$ and contain more information
than the ramification indices.

\begin{theorem}\label{thm:kummer-normalizer}
The deck group $N$ is central in $M$. There is an exact sequence
\[
 1\longrightarrow C_n\longrightarrow M\longrightarrow L_{n,m}
 \longrightarrow1,
\]
where
\[
 L_{n,m}=
 \{\tau\in\operatorname{PGL}_2(k):w(\tau(P))=w(P)
                       \text{ for every }P\in\mathbb P^1(k)\}.
\]
More explicitly, if
\[
 H_m=\operatorname{Stab}_{\operatorname{PGL}_2(k)}(\mathcal Z_m),\qquad
 K_m=\operatorname{Stab}_{\operatorname{PGL}_2(k)}
                  (\mathcal Z_m\cup\{\infty\}),
\]
then
\[
 L_{n,m}=
 \begin{cases}
 H_m,&n\mid m,\\
 K_m,&n\mid(m+1),\\
 \langle x\mapsto-x\rangle\cong C_2,
       &n\nmid m\text{ and }n\nmid(m+1).
 \end{cases}
\]
\end{theorem}
\begin{proof}
An element of $M$ induces a fractional linear transformation $\tau$ on
$F^N=k(x)$. The eigenspaces of $N$ in $F$ are the one-dimensional
$k(x)$-spaces $k(x)y^j$. Conjugation by an element of $M$ permutes
them, so its action has the form
\[
 x\longmapsto\tau(x),\qquad y\longmapsto u(x)y^a,
 \quad u(x)\in k(x)^*,\quad a\in(\mathbb Z/n\mathbb Z)^*.
\]
Hence
\begin{equation}\label{eq:kummer-lifting}
 \varphi_m(\tau(x))=u(x)^n\varphi_m(x)^a.
\end{equation}
Every finite zero must map to a finite zero or to infinity: the valuation
of the right side at such a point is nonzero modulo $n$. Since $m\geq2$,
at least one finite zero maps to another finite zero. Taking valuations at
that point gives $a\equiv1\pmod n$. Thus every element of $M$ commutes with
$\nu$, and its induced transformation preserves $w$.

Conversely, if $\tau$ preserves $w$, the divisor of
$\varphi_m(\tau(x))/\varphi_m(x)$ is divisible by $n$. Since every divisor
of degree zero on $\mathbb P^1$ is principal and every element of $k^*$ is
an $n$th power, this quotient equals $u(x)^n$ for some $u\in k(x)^*$.
The rule $(x,y)\mapsto(\tau(x),u(x)y)$ is then a field automorphism, hence
an automorphism of the smooth projective curve. There are exactly $n$
such lifts. This proves the exact sequence.

If $n\mid m$, infinity has weight zero and $w$ distinguishes precisely
$\mathcal Z_m$. If $n\mid(m+1)$, infinity and all finite zeros have weight one.
In the remaining case, infinity is the unique point with weight $-m$,
which is neither zero nor one. It is therefore fixed, and $\tau=ax+b$ is
affine and preserves $\mathcal Z_m$. Comparing coefficients in
\[
 \varphi_m(ax+b)=a^m\varphi_m(x)
\]
uses $\varphi_m(x)=x^m-mx^{m-2}+\cdots$. The coefficient of $x^{m-1}$
gives $b=0$, because $p\nmid m$, and the coefficient of $x^{m-2}$ then
gives $a^2=1$. Both possibilities occur since
$\varphi_m(-x)=(-1)^m\varphi_m(x)$.
\end{proof}

\begin{remark}
Suppose $\gcd(n,m)=1$ and $n\nmid(m+1)$. All branch points are then
totally ramified, but infinity has a different Kummer weight. Every
transformation induced by $M$ must therefore fix infinity.
\end{remark}

\begin{proposition}\label{prop:kummer-affine}
The inverse image in $M$ of $\langle x\mapsto-x\rangle$ is
\[
 \mathcal A_{n,m}=
 \{(x,y)\mapsto(\epsilon x,cy):\epsilon\in\{1,-1\},\
                         c^n=\epsilon^m\}.
\]
It has order $2n$, and
\[
 \mathcal A_{n,m}\cong
 \begin{cases}
 C_n\times C_2,&m\text{ even},\\
 C_{2n},&m\text{ odd}.
 \end{cases}
\]
In the odd case one generator is $(x,y)\mapsto(-x,\zeta_{2n}y)$.
When $n$ is odd these two abstract groups are isomorphic.
\end{proposition}
\begin{proof}
The formula follows from $\varphi_m(-x)=(-1)^m\varphi_m(x)$ and the
lifting equation. If $m$ is even, the involution $(x,y)\mapsto(-x,y)$ is
independent of $N$. If $m$ is odd, the displayed generator has order
$2n$ and its square generates $N$.
\end{proof}

\section{Explicit exceptional normalizers}

For $m\geq6$, equation~\eqref{eq:hyp-companion-input} and
Proposition~\ref{prop:hyp-odd-stabilizer} give the following two
stabilizer classifications:
\begin{align*}
 H_m&\cong
 \begin{cases}
 D_m,&\operatorname{char}(k)=p>0,\ 2m=p^r+1,\\
 S_3,&\operatorname{char}(k)=p>0,\ 4m=p^r+1,\\
 C_2,&\text{otherwise};
 \end{cases}\\
 K_m&\cong
 \begin{cases}
 D_{m+1},&\operatorname{char}(k)=p>0,\ 2m+1=p^r,\\
 C_2,&\text{otherwise}.
 \end{cases}
\end{align*}
In every $C_2$ row the group is generated by $x\mapsto-x$.

For $n\mid d$, define a group
\begin{equation}\label{eq:End-presentation}
 \mathcal E_{n,d}=\left\langle\nu,R,S\ \middle|\
 \nu^n=R^d=S^2=1,\ [\nu,R]=[\nu,S]=1,
 \ SRS=\nu^{-1}R^{-1}\right\rangle.
\end{equation}

\begin{proposition}\label{prop:End}
The group $\mathcal E_{n,d}$ has order $2nd$ and is isomorphic to
$(C_n\times C_d)\rtimes C_2$, where the involution fixes $\nu$ and sends
$R$ to $\nu^{-1}R^{-1}$. It is a central extension
\[
 1\longrightarrow\langle\nu\rangle\longrightarrow \mathcal E_{n,d}
   \longrightarrow D_d\longrightarrow1.
\]
This central extension splits if and only if $n$ is odd. For odd $n$,
$\mathcal E_{n,d}\cong C_n\times D_d$. For $n=2$ and $d\geq6$, the group $\mathcal E_{2,d}$
is $\mathcal G_d$ from \eqref{eq:hyp-GN}.
\end{proposition}
\begin{proof}
On the smooth model of $v^n=u^d+1$, choose a primitive $d$th root
$\zeta_d$ and put $\zeta_n=\zeta_d^{d/n}$. The formulas
\[
 \nu(u,v)=(u,\zeta_n v),\qquad
 R(u,v)=(\zeta_d u,v),\qquad
 S(u,v)=\left(u^{-1},\frac{v}{u^{d/n}}\right)
\]
give automorphisms and satisfy \eqref{eq:End-presentation}. The elements
$\nu^aR^bS^c$, with $0\leq a<n$, $0\leq b<d$ and $c\in\{0,1\}$,
are distinct. Conversely, the relations reduce every word to this form,
so the presentation has exactly $2nd$ elements.

A complement to $\langle\nu\rangle$ would contain lifts
$R'=\nu^aR$ and $S'=\nu^bS$ of the standard rotation and reflection in
$D_d$. The equality $S'R'S'^{-1}=R'^{-1}$ forces
$2a\equiv1\pmod n$. Thus no complement exists for even $n$. If $n$ is
odd, choose $a$ with $2a\equiv1\pmod n$. Then $R'=\nu^aR$ has order
$d$, because $n\mid d$, and $\langle R',S\rangle\cong D_d$ is a
complement. Since $\nu$ is central, the resulting product is direct.
\end{proof}

\begin{theorem}\label{thm:explicit-normalizers}
Assume $m\geq6$. Then
\[
 M\cong
 \begin{cases}
 \mathcal E_{n,m},&n\mid m,\quad 2m=p^r+1,\\
 C_n\times S_3,&n\mid m,\quad 4m=p^r+1,\\
 \mathcal E_{n,m+1},&n\mid(m+1),\quad 2m+1=p^r,\\
 \mathcal A_{n,m},&\text{otherwise}.
 \end{cases}
\]
The first three rows require positive characteristic and $r\geq1$.
Their orders are respectively $2nm$, $6n$, and $2n(m+1)$.
\end{theorem}
\begin{proof}
Theorem~\ref{thm:kummer-normalizer}, the stabilizer classifications
above and Proposition~\ref{prop:kummer-affine} give the last row. It
remains to identify the groups in the three exceptional cases. Put
\[
 z=\frac{x-2}{x+2},\qquad
 P_m(z)=\sum_{j=0}^m\binom{2m}{2j}z^j.
\]
The defining identity for Chebyshev polynomials gives
\begin{equation}\label{eq:Cheb-binary}
 (1-z)^m\varphi_m\!\left(\frac{2(1+z)}{1-z}\right)=2P_m(z).
\end{equation}
If $2m=q+1$ with $q=p^r$, the Frobenius identity gives
$P_m(z)=1+z^m$. Choosing $c^n=2$ and setting
$v=c^{-1}(1-z)^{m/n}y$ gives $v^n=1+z^m$. The lifts in
Proposition~\ref{prop:End} have order $2nm$, which is the order prescribed
by the exact sequence. This proves the first row.

If $2m+1=q=p^r$, the congruences
$\binom{q-1}{j}=(-1)^j$ in $k$ give
\[
 P_m(z)=1+z+\cdots+z^m=\frac{1-z^{m+1}}{1-z}.
\]
Thus $v=c^{-1}(1-z)^{(m+1)/n}y$ satisfies
$v^n=1-z^{m+1}$. Scaling $z$ changes this to the model in
Proposition~\ref{prop:End}, giving the third row.

For the second row put $t=m/n$ and
\[
 h(x)=\frac{2x+12}{2-x}.
\]
In the exceptional $S_3$ stabilizer case the standard generators are
$x\mapsto-x$ and $h$. They satisfy $h^3=1$ and
$(-x)\circ h\circ(-x)=h^{-1}$. Their lifts
\begin{equation}\label{eq:S3-Kummer-lifts}
 s(x,y)=(-x,(-1)^m y),\qquad
 T(x,y)=\left(h(x),\frac{(-4)^t y}{(2-x)^t}\right)
\end{equation}
preserve the curve and satisfy $s^2=T^3=1$ and $sTs=T^{-1}$.
Here the polynomial identity needed for $T$ is
\[
 (2-x)^m\varphi_m(h(x))=(-4)^m\varphi_m(x).
\]
Indeed, preservation of the zero set makes the two degree-$m$ polynomials
proportional. Comparing their leading coefficients gives the scalar $2$.
Also $q=4m-1\equiv3\pmod4$, and if $i^2=-1$ then
$(1+i)^{q+1}=2$ and $(1+i)^4=-4$, so $(-4)^m=2$.
The matrix
\[
 B=\begin{pmatrix}2&12\\-1&2\end{pmatrix}
\]
satisfies $B^3=-64I$ and
$\operatorname{diag}(-1,1)B\operatorname{diag}(-1,1)=-\tfrac14B^2$.
These equalities verify the relations for the lifts in
\eqref{eq:S3-Kummer-lifts}. The subgroup they generate maps isomorphically
onto $S_3$, and hence meets the central deck group trivially. It gives
the asserted direct product of order $6n$.
\end{proof}

\section{When the normalizer is the full automorphism group}

\begin{theorem}\label{thm:strong-branching-Cheb}
If $m>2n$, then the deck group $N$ is normal in
$\operatorname{Aut}(X_{n,m})$. Consequently $M=\operatorname{Aut}(X_{n,m})$.
This assertion is valid for arbitrary $n$, including composite $n$, and
in every characteristic allowed above.
\end{theorem}
\begin{proof}
Let $\sigma\in\operatorname{Aut}(F)$ and put
$N'=\sigma N\sigma^{-1}$ and $I=N\cap N'$. Write $|I|=e$ and $a=n/e$.
The field
\[
 L=F^NF^{N'}=k(x,\sigma(x))
\]
is an intermediate field of the Galois extension $F/F^N$.
The extension $F/L$ is therefore Galois, with
\[
 \operatorname{Gal}(F/L)=N\cap N'=I.
\]
It follows that $L=F^I=k(x,y^e)$ and that
\[
 [L:k(x)]=[L:k(\sigma(x))]=a.
\]
If $N'\neq N$, then $a\geq2$. The two rational functions $x$ and
$\sigma(x)$ generate $L$, and hence give a birational map of its smooth
curve to an integral curve of bidegree $(a,a)$ in
$\mathbb P^1\times\mathbb P^1$. Its arithmetic genus is $(a-1)^2$, so
(equivalently, by Castelnuovo's inequality \cite[Theorem~3.11.3]{Stichtenoth})
\begin{equation}\label{eq:CS-Cheb}
 g(L)\leq(a-1)^2.
\end{equation}
On the other hand $L=k(x,y^e)$ has equation
$(y^e)^a=\varphi_m(x)$, and Proposition~\ref{prop:kummer-genus} gives
\[
 g(L)=\frac{(a-1)(m-1)-\gcd(a,m)+1}{2}
 \geq\frac{(a-1)(m-2)}{2}
 >(a-1)^2,
\]
since $m>2n\geq2a$. This contradicts \eqref{eq:CS-Cheb}, so $N'=N$.
\end{proof}

\begin{remark}
Kontogeorgis~\cite[Proposition~1]{Kontogeorgis1999} proved the normality
criterion $s>2n$ for a cyclic degree-$n$ cover of the rational function
field fully ramified at $s$ places. See also
\cite[Proposition~2.1]{AntoniadisKontogeorgis} (arXiv numbering).
We included the proof to state explicitly how the Castelnuovo--Severi
inequality applies for composite $n$ and in positive characteristic.
\end{remark}

\begin{corollary}\label{cor:full-Cheb-high-m}
If $m\geq6$ and $m>2n$, Theorem~\ref{thm:explicit-normalizers} is a
classification of the full group $\operatorname{Aut}(X_{n,m})$.
In characteristic zero this full group is $C_n\times C_2$ for even $m$
and $C_{2n}$ for odd $m$.
\end{corollary}

\begin{proposition}\label{prop:large-n-preserves-y}
If $n>2m$, every automorphism of $X_{n,m}$ preserves the rational
subfield $k(y)$.
\end{proposition}
\begin{proof}
Introduce a field
\[
 E=k(t,y),\qquad y^n=t^m+t^{-m},\qquad x=t+t^{-1}.
\]
The rational function $t^m+t^{-m}$ has simple zeros, so its Kummer
extension has degree $n$. The transformations
\[
 r(t,y)=(\zeta_m t,y),\qquad s(t,y)=(t^{-1},y)
\]
generate a dihedral group of order $2m$. The polynomial
$t^{2m}-y^nt^m+1$ has degree $2m$ in $t$, so
$E/k(y)$ is Galois with exactly this group. We have
$F=k(x,y)=E^{\langle s\rangle}$.

Every subgroup of $D_m=\langle r,s\rangle$ containing $s$ has the form
$J_b=\langle r^b,s\rangle$ for some divisor $b$ of $m$. Its fixed field is
\[
 E^{J_b}
 =k\bigl(y,t^{m/b}+t^{-m/b}\bigr)
 =k\bigl(y,\varphi_{m/b}(x)\bigr).
\]
Writing $u=\varphi_{m/b}(x)$, the Chebyshev composition identity gives
$y^n=\varphi_b(u)$. Thus this fixed field is the function field of
$X_{n,b}$ and has degree $b$ over $k(y)$.

Now let $\sigma\in\operatorname{Aut}(F)$ and put
$L=k(y,\sigma(y))$. It is an intermediate field of $F/k(y)$ and hence
equals $E^{J_b}$ for a divisor $b$ of $m$. Moreover
\[
 [L:k(y)]=[L:k(\sigma(y))]=b,
\]
because both $k(y)$ and $k(\sigma(y))$ have degree $m$ in $F$. If these
two rational subfields differ, then $b\geq2$. Their generators give a
birational image of the curve of $L$ of bidegree $(b,b)$, so its genus
is at most $(b-1)^2$. But
\[
 g(L)=\frac{(n-1)(b-1)-\gcd(n,b)+1}{2}
 \geq\frac{(n-2)(b-1)}2
 >(b-1)^2,
\]
as $n>2m\geq2b$. This contradiction proves $k(\sigma(y))=k(y)$.
\end{proof}

\begin{theorem}\label{thm:full-Cheb-high-n}
Assume $n>2m$ and $m/\gcd(m,n)>2$. Then
\[
 \operatorname{Aut}(X_{n,m})=\mathcal A_{n,m},
\]
so the full group has order $2n$. It is $C_n\times C_2$ for even $m$
and $C_{2n}$ for odd $m$.
\end{theorem}
\begin{proof}
By Proposition~\ref{prop:large-n-preserves-y}, every automorphism induces
a fractional linear transformation on the $y$-line. The finite branch
values of the degree-$m$ map $y:X_{n,m}\to\mathbb P^1$ are
\[
 W=\{a\in k:a^{2n}=4\}.
\]
The assumption $m/\delta>2$ implies $m\geq3$, so both Chebyshev critical
values $2$ and $-2$ occur. All ramification indices above $W$ are one or
two. The pole orders of $x$ and $y$ at each of the $\delta$ points above
infinity are $n/\delta$ and $m/\delta$, respectively. Thus the
ramification indices of the map to the $y$-line above infinity are
all $m/\delta>2$. Hence infinity is fixed
by every induced transformation, which is therefore $y\mapsto ay+b$.
Comparison of the coefficient of $y^{2n-1}$ in
\[
 (ay+b)^{2n}-4=a^{2n}(y^{2n}-4)
\]
gives $b=0$, since $p\nmid2n$, and comparison of constants gives
$a^{2n}=1$.

If $m$ is even, each fiber with $y^n=2$ has exactly two unramified
points, while each fiber with $y^n=-2$ has none. This follows from
\[
 \varphi_m(x)-2=(x^2-4)Q(x)^2,\qquad
 \varphi_m(x)+2=R(x)^2
\]
for suitable polynomials $Q,R$ whose zeros are simple and avoid the
unramified points. These assertions follow either from the substitution
$x=t+t^{-1}$ or from the identities in
\eqref{eq:hyp-cheb-identities}; the assumptions $p\nmid2m$ ensure
that each critical point has ramification index exactly two. The induced transformation must preserve the two
types, so $a^n=1$. The image on the $y$-line therefore has order at most
$n$. For odd $m$, the analogous factorizations are
$\varphi_m(x)-2=(x-2)Q(x)^2$ and
$\varphi_m(x)+2=(x+2)R(x)^2$. Thus the two critical fibers each have
one unramified point, and the bound on the image is $2n$.

The kernel is $\operatorname{Aut}(F/k(y))$. With the dihedral Galois
extension used above it is
\[
 \operatorname{N}_{D_m}(\langle s\rangle)/\langle s\rangle.
\]
This group has order two when $m$ is even and is trivial when $m$ is
odd. Thus in either parity the full automorphism group has order at most
$2n$. It already contains $\mathcal A_{n,m}$, which has order $2n$, proving the
assertion.
\end{proof}

\begin{lemma}\label{lem:Kummer-no-inversion}
Suppose $m\geq3$ and $m\mid n$. An automorphism of $X_{n,m}$ that
preserves $k(y)$ cannot induce $y\mapsto a/y$ on the $y$-line.
\end{lemma}
\begin{proof}
The branch set of the map to the $y$-line is
$W=\{y:y^{2n}=4\}$. An inversion preserving $W$ has
$a^n=4$ or $-4$.
Use the dihedral extension $E$ of
Proposition~\ref{prop:large-n-preserves-y}. Since $m\geq3$, the
reflection subgroup $\langle s\rangle$ has trivial core in $D_m$, so
$E$ is the Galois closure of $F/k(y)$. An automorphism of $F$ preserving
$k(y)$ extends to $E$ and normalizes $D_m$. The rotation subgroup
$\langle r\rangle\cong C_m$ is characteristic in $D_m$, since an
element of order $m>2$ cannot be a reflection. Thus the automorphism
also preserves
\[
 F_0=E^{\langle r\rangle}=k(u,y),\qquad
 u=t^m,\qquad y^n=u+u^{-1}.
\]
If it induces $y\mapsto a/y$, its image $u'$ satisfies
\[
 u'+u'^{-1}=\frac{a^n}{u+u^{-1}}.
\]
For $i^2=-1$, the two solutions in $k(u)$ are
$i(u-i)/(u+i)$ and its inverse when $a^n=4$; their negatives are the
solutions when $a^n=-4$. Consequently the valuations of $u'$ at
$u=i$ and $u=-i$ are $1,-1$ in one order or the other.

An extension to $E=F_0(t)$, $t^m=u$, that normalizes its deck $C_m$
must have the form $t\mapsto b t^c$, where $b\in F_0^*$ and
$c\in(\mathbb Z/m\mathbb Z)^*$. Hence
\[
 b^m=\frac{u'}{u^c}\in k(u).
\]
The extension $F_0/k(u)$ is cyclic of degree $n$, with equation
$y^n=(u^2+1)/u$. If $b^m\in k(u)$, applying a generator of this cyclic
group shows that $b$ is an eigenvector. It follows that
$b=h(u)y^e$ for some $h\in k(u)^*$ and $0\leq e<n$. Since $b^m$ is
fixed, $me=nj$ for an integer $j$, and therefore
\[
 \frac{u'}{u^c}
 =h(u)^m\left(\frac{u^2+1}{u}\right)^j.
\]
Taking valuations at $u=i$ and $u=-i$ gives simultaneously
$j\equiv1\pmod m$ and $j\equiv-1\pmod m$. Thus $m\mid2$, contrary
to $m\geq3$. Thus no inversion of the $y$-line lifts.
\end{proof}

\begin{lemma}\label{lem:roots-unity-stabilizer}
Let $d\geq3$ and assume $\operatorname{char}(k)\nmid d$.
In characteristic zero, or in characteristic $p>0$ if $d\ne p^r+1$
for every $r\geq1$, one has
\[
 \operatorname{Stab}_{\operatorname{PGL}_2(k)}(\mu_d)
 =\langle z\mapsto\zeta_dz,\ z\mapsto z^{-1}\rangle\cong D_d.
\]
If $d=q+1$ with $q=p^r$, this stabilizer is conjugate to
$\operatorname{PGL}_2(\mathbb F_q)$.
\end{lemma}
\begin{proof}
A fractional linear transformation $(az+b)/(cz+e)$ preserving $\mu_d$
satisfies
\begin{equation}\label{eq:root-unity-binary}
 (az+b)^d-(cz+e)^d=\lambda(z^d-1),\qquad \lambda\ne0.
\end{equation}
Since $d\ne0$ in $k$, comparison of the coefficients of $z$ and
$z^{d-1}$ shows that if any of $a,b,c,e$ vanishes, the opposite entry
vanishes as well. Such a transformation is a scaling or an inversion,
and \eqref{eq:root-unity-binary} puts it in the displayed dihedral group.

Suppose $abce\ne0$ and put $A=a/c$ and $B=b/e$. Each $j$ with
$1\leq j\leq d-1$ and $\binom dj\ne0$ in $k$ satisfies
\[
 A^jB^{d-j}=1.
\]
Two consecutive such integers $j$ force $A=B$, contradicting
$ae-bc\ne0$. In characteristic zero, $j=1,2$ supply the contradiction.
In characteristic $p$, write $d=\sum d_rp^r$ in base $p$. The expansion
\[
 (1+T)^d=\prod_r(1+T^{p^r})^{d_r}
\]
shows that $\binom dj$ is nonzero precisely when every base-$p$ digit
of $j$ is at most the corresponding digit of $d$. Since $p\nmid d$,
$d_0\geq1$. If $d_0\geq2$, the indices $1,2$ give a contradiction,
so $d_0=1$. Choose $r\geq1$ with $d_r>0$. Unless $d=p^r+1$, the two
internal indices $p^r,p^r+1$ have nonzero coefficients, again a
contradiction. Thus a nonmonomial transformation requires $d=q+1$.

For $d=q+1$, choose $\theta\in\mathbb F_{q^2}\setminus\mathbb F_q$.
The transformation
\[
 z\longmapsto\frac{z-\theta}{z-\theta^q}
\]
maps $\mathbb P^1(\mathbb F_q)$ bijectively to $\mu_{q+1}$. The
stabilizer of $\mathbb P^1(\mathbb F_q)$ is exactly
$\operatorname{PGL}_2(\mathbb F_q)$. Indeed, a projective transformation
is determined by its values at three distinct points. Any two ordered
triples of distinct $\mathbb F_q$-rational points are related by a
unique element of $\operatorname{PGL}_2(\mathbb F_q)$.
\end{proof}

We use the following consequence of Dickson's subgroup classification.

\begin{lemma}\label{lem:Dickson}
Let $q$ be a power of an odd prime with $q\geq5$, and let
$H\leq\operatorname{PGL}_2(\mathbb F_q)$ be a subgroup containing a cyclic
subgroup $T$ of order $q+1$. Then either
$H=\operatorname{PGL}_2(\mathbb F_q)$, or $H$ normalizes $T$.
\end{lemma}
\begin{proof}
Every element of $\operatorname{PGL}_2(\mathbb F_q)$ has order $p$ or
order dividing $q-1$ or $q+1$. Thus $q+1$ is the largest possible
element order. Apply Dickson's classification \cite{Dickson1901};
see also Valentini and Madan~\cite{ValentiniMadan}. The possibilities
are cyclic groups, dihedral groups, elementary abelian $p$-groups and
their affine extensions, the groups $A_4,S_4,A_5$, and the groups
$\operatorname{PSL}_2(\mathbb F_{q'})$ and
$\operatorname{PGL}_2(\mathbb F_{q'})$ over subfields of $\mathbb F_q$.

An affine subgroup has no element of order $q+1$. Indeed, its
prime-to-$p$ element orders divide $q-1$. The groups $A_4,S_4,A_5$
have no element of order at least six. An element of
$\operatorname{PSL}_2(\mathbb F_{q'})$ has order $p$ or order dividing
$(q'\pm1)/2$. If $q'<q$, every element of
$\operatorname{PGL}_2(\mathbb F_{q'})$ has order at most $q'+1<q+1$.
These cases are therefore excluded.

It follows that $H$ is cyclic, dihedral, or the full group.
If $H$ is cyclic, the bound on element orders gives $H=T$.
If $H$ is dihedral, a generator of $T$ cannot be a reflection, since
it has order $q+1\geq6$. Thus $T$ lies in the cyclic subgroup of
index two. The same bound on element orders shows that this cyclic
subgroup is $T$. Hence $H$ normalizes $T$.
\end{proof}

\begin{theorem}\label{thm:full-Cheb-high-n-sharpened}
Assume $m\geq3$ and $n>2m$, and let $k$ be algebraically closed of
characteristic zero or of characteristic prime to $2nm$. Then
\[
 \operatorname{Aut}(X_{n,m})=\mathcal A_{n,m}.
\]
In particular the full group has order $2n$; it is $C_n\times C_2$ for
even $m$ and $C_{2n}$ for odd $m$.
\end{theorem}
\begin{proof}
Put $\delta=\gcd(m,n)$. Theorem~\ref{thm:full-Cheb-high-n} handles
$m/\delta>2$, leaving $m/\delta=1$ or $2$. The projection to $y$ is
preserved in every case by Proposition~\ref{prop:large-n-preserves-y}.

First suppose $m/\delta=1$. The branch set is $W=\{y:y^{2n}=4\}$.
We separate the case where $m$ is odd and $2n=p^r+1$.
Outside that case, the argument is as follows.
If $2n$ is not one more than a power of the characteristic,
Lemma~\ref{lem:roots-unity-stabilizer} shows that every induced
transformation is $y\mapsto ay$ or $y\mapsto a/y$. If instead
$2n=p^r+1$ and $m$ is even, the critical fiber profiles still distinguish
$B_+=\{y:y^n=2\}$. Now $n\equiv1/2\pmod p$, so $n$ cannot equal
$p^s+1$ for any $s\geq1$. Applying the same lemma to $B_+$ therefore
gives the same two types of transformations. Lemma~\ref{lem:Kummer-no-inversion} excludes the inverse type.
Only rotations remain, and the image and kernel bounds give order $2n$.

Now suppose $m/\delta=1$, $m$ is odd, and $2n=q+1$ with $q=p^r$.
Since $m\mid n$ and $n>2m$, we have $n\geq3m\geq9$. Thus $q\geq17$.
By Lemma~\ref{lem:roots-unity-stabilizer}, the stabilizer
$\operatorname{Stab}(W)$ is conjugate to $\operatorname{PGL}_2(\mathbb F_q)$.
Let $\Gamma\leq\operatorname{Stab}(W)$ be the group of transformations of
the $y$-line induced by $\operatorname{Aut}(X_{n,m})$. Since $m$ is odd,
the kernel calculation in the proof of
Theorem~\ref{thm:full-Cheb-high-n} gives
$\operatorname{Aut}(F/k(y))=1$. Hence
$\operatorname{Aut}(X_{n,m})\cong\Gamma$. The generator
$(x,y)\mapsto(-x,\zeta_{2n}y)$ of $\mathcal A_{n,m}$ induces
$y\mapsto\zeta_{2n}y$, so $\Gamma$ contains the cyclic group
\[
 T=\{y\mapsto\zeta y:\zeta^{2n}=1\}
\]
of order $2n=q+1$. Choose $a\in k$ with $a^n=4$. The inversion
$y\mapsto a/y$ preserves $W$, but it does not lie in $\Gamma$ by
Lemma~\ref{lem:Kummer-no-inversion}. Hence
$\Gamma\ne\operatorname{Stab}(W)$, and Lemma~\ref{lem:Dickson}, applied
after conjugating $\operatorname{Stab}(W)$ to
$\operatorname{PGL}_2(\mathbb F_q)$, shows that $\Gamma$ normalizes $T$.
Every element of the normalizer of $T$ in $\operatorname{PGL}_2(k)$
permutes the two fixed points $0,\infty$ of $T$. If it fixes both, it is
a scaling $y\mapsto by$ preserving $W$, so $b^{2n}=1$ and it lies in $T$.
If it interchanges them, it is an inversion, which is excluded by
Lemma~\ref{lem:Kummer-no-inversion}. Therefore $\Gamma=T$ and
$|\operatorname{Aut}(X_{n,m})|=2n=|\mathcal A_{n,m}|$.

Finally, suppose $m/\delta=2$. Then $m$ is even. The fiber above infinity consists
of $\delta=m/2$ points, all with ramification index two, just as
each fiber with $y^n=-2$. By contrast, each fiber with $y^n=2$
has two unramified points. The critical fiber profiles therefore distinguish
\[
 B_+=\{y:y^n=2\},\qquad
 S=\{y:y^n=-2\}\cup\{\infty\}.
\]
Every induced transformation preserves both sets. Unless
$n=q+1$ with $q=p^r$, Lemma~\ref{lem:roots-unity-stabilizer} applied to
$B_+$ makes it a scaling or an inversion. An inversion sends
$\infty\in S$ to $0\notin S$, so it cannot occur.

It remains to consider $n=q+1$. Write the induced transformation as
$\gamma(y)=(ay+b)/(cy+e)$. Preservation of the two homogeneous branch forms
gives
\begin{align*}
 (ay+b)^n-2(cy+e)^n&=\lambda(y^n-2),\\
 (cy+e)\bigl((ay+b)^n+2(cy+e)^n\bigr)&=\mu(y^n+2),
\end{align*}
with $\lambda,\mu\ne0$. Eliminating $(ay+b)^n$ gives
\[
 \lambda(cy+e)(y^{q+1}-2)+4(cy+e)^{q+2}
 =\mu(y^{q+1}+2).
\]
Since
\[
 (cy+e)^{q+2}=(c^qy^q+e^q)(cy+e)^2,
\]
its coefficient of $y^q$ is $4c^qe^2$. Here $q\geq3$, so that
power is distinct from the other relevant powers. If $c\ne0$, then
$e=0$; the determinant condition gives $b\ne0$. The coefficient of
$y^q$ in the first branch-form identity now gives $a^qb=0$, hence
$a=0$. This would again be an inversion, contradicting preservation of
$S$. Therefore $c=0$. The transformation is affine and preserves
$B_+$, so comparison of its degree-$n-1$ coefficient makes it a scaling.
The image and kernel bounds in Theorem~\ref{thm:full-Cheb-high-n}
finish the proof.
\end{proof}

\begin{remark}\label{rem:high-n-Dickson}
The branch stabilizer can be much larger than the group that lifts.
For $(n,m,p)=(9,3,17)$, it is conjugate to
$\operatorname{PGL}_2(\mathbb F_{17})$, of order $4896$.
The theorem gives $\operatorname{Aut}(X_{9,3})\cong C_{18}$.
For $(n,m,p)=(15,5,29)$, the corresponding orders are $24360$ and $30$.
In both cases, the transformations that lift to $X_{n,m}$ are exactly
the scalings in $T$. This follows from
Lemmas~\ref{lem:Kummer-no-inversion} and~\ref{lem:Dickson}.
\end{remark}

\subsection{The quadratic Chebyshev polynomial}
\label{sec:quadratic-Chebyshev}

The case $m=2$ of a Chebyshev--Kummer curve is a binomial hyperelliptic
curve after interchanging the roles of the coordinates. We determine
its full group, including the cases where its distinguished cyclic
deck group is not normal. Throughout this subsection $k$ is algebraically
closed of characteristic zero or of odd characteristic $p$.

In positive characteristic $p$, for a power $q=p^e\geq5$, define the group
\begin{equation}\label{eq:binomial-Rq-group}
 \mathcal R_q=
 \frac{\{(A,\eta):A\in\operatorname{GL}_2(\mathbb F_q),\
                 \eta\in k^*,\ \eta^2=\det A\}}
      {\{(\lambda I,\lambda^{(q+1)/2}):
                         \lambda\in\mathbb F_q^*\}},
\end{equation}
where multiplication of pairs is $(A,\eta)(B,\xi)=(AB,\eta\xi)$.
The denominator is a central subgroup. This group will occur as the
full automorphism group of the curve $R_q:V^2=U^q-U$.

\begin{theorem}\label{thm:binomial-hyperelliptic}
Let $N\geq5$, assume $p\nmid N$ in positive characteristic, and let
$Y_N$ be the smooth projective curve with equation
\[
 v^2=u^N-1.
\]
Then
\[
 \operatorname{Aut}(Y_N)\cong
 \begin{cases}
  C_{2N},&N\text{ odd},\\
  \mathcal G_N,&N\text{ even and }N\ne p^e+1
                              \text{ for every }e\geq1,\\
  \mathcal R_q,&N=q+1,\quad q=p^e.
 \end{cases}
\]
Here $\mathcal G_N$ is the group in \eqref{eq:hyp-GN}, of order $4N$.
The third case requires positive characteristic; its group has order
$2q(q^2-1)$ and fits into a central exact sequence
\[
 1\longrightarrow C_2\longrightarrow\mathcal R_q
 \longrightarrow\operatorname{PGL}_2(\mathbb F_q)
 \longrightarrow1.
\]
In characteristic zero, only the first two cases occur.
\end{theorem}

For the classification in characteristic zero, see Brandt and
Stichtenoth~\cite{BrandtStichtenoth}. The proof below also covers odd
positive characteristic.

\begin{proof}
The polynomial $u^N-1$ is separable. The genus is
$\lfloor(N-1)/2\rfloor\geq2$, so the full group is a central extension
by the hyperelliptic involution of the stabilizer of its branch set.
That branch set is $\mu_N\cup\{\infty\}$ when $N$ is odd and
$\mu_N$ when $N$ is even.

\emph{Odd $N$.}
Put $\Omega=\mu_N\cup\{\infty\}$ and
$K=\operatorname{Stab}_{\operatorname{PGL}_2(k)}(\Omega)$.
An affine map $u\mapsto au+b$ in $K$ satisfies
$(au+b)^N-1=a^N(u^N-1)$. Its coefficients of degrees $N-1$ and zero
give $b=0$ and $a^N=1$. Thus $K_\infty$ is the rotation group $C_N$.
We show that no nonaffine element exists.

Suppose $M\in K$ moves infinity. Its image of infinity and its pole
both belong to $\mu_N$. Composing on the left and right by suitable
rotations independently makes both points equal to $1$. Hence we may
write
\[
 M(u)=1+\frac{s}{u-1},\qquad s\ne0.
\]
The homogeneous form defining $\Omega$ has degree $N+1$. Therefore
\begin{equation}\label{eq:binomial-odd-identity}
 (u-1)^{N+1}\bigl(M(u)^N-1\bigr)=\lambda(u^N-1),
 \qquad \lambda\in k^*.
\end{equation}
At each $\beta\in\mu_N\setminus\{1\}$, differentiate this identity
twice and divide its second derivative by its first derivative.
Using $M''/M'=-2/(u-1)$ gives
\[
 \frac{N-1}{\beta}
 =\frac{2N}{\beta-1}
   +M'(\beta)\frac{N-1}{M(\beta)}.
\]
All denominators are nonzero: $\beta$ and $M(\beta)$ are nonzero
finite elements of $\Omega$. The numerator of the difference between
the two sides, as rational functions of $u$ with common denominator
$u(u-1)(u+s-1)$, is
\begin{equation}\label{eq:binomial-quadratic-obstruction}
 -(N+1)u^2+2(1-s)u+(N-1)(1-s).
\end{equation}
It vanishes at the $N-1\geq4$ distinct points just considered.
Consequently it is zero, so $N+1=0$ in $k$ and $s=1$.
The first equality is already impossible in characteristic zero.

In positive characteristic, the transformation is now
$M(u)=u/(u-1)$. In the coordinate $w=1/u$, it becomes
$w\mapsto1-w$, and $\Omega$ becomes $\{0\}\cup\mu_N$, the simple
zero set of
\[
 F(w)=w^{N+1}-w.
\]
Since $N$ is odd, both $F(1-w)$ and $F(w)$ are monic. Preservation of
their zero sets would give $F(1-w)=F(w)$. The coefficient of $w$ on
the left is $-N=1$, because $N+1=0$ in $k$, whereas that on the
right is $-1$. This contradicts the assumption that the
characteristic is odd. We have proved $K=C_N$ in every allowed
characteristic.

The rotation $(u,v)\mapsto(\zeta_Nu,v)$ and the hyperelliptic
involution $(u,v)\mapsto(u,-v)$ commute and generate a group of
order $2N$. It is the full group, and it is cyclic because $N$ is odd.

\emph{Even $N$, outside the Frobenius case.}
Lemma~\ref{lem:roots-unity-stabilizer} gives
$\operatorname{Stab}(\mu_N)=D_N$ unless $N=p^e+1$.
Over $k$, the curve $v^2=u^N-1$ is isomorphic, by rescaling the
coordinates, to $w^2=z^N+1$. The transformations in
\eqref{eq:hyp-GN-actions} form a subgroup $\mathcal G_N$ of order
$4N$. This is the full order supplied by the hyperelliptic exact
sequence, proving the second case.

\emph{Even $N=q+1$, where $q=p^e$.}
Choose $\theta\in\mathbb F_{q^2}\setminus\mathbb F_q$ and choose
$c_0\in k$ with $c_0^2=\theta^q-\theta$. The change of variables
\begin{equation}\label{eq:binomial-Rq-change}
 u=\frac{U-\theta}{U-\theta^q},\qquad
 v=\frac{c_0V}{(U-\theta^q)^{(q+1)/2}}
\end{equation}
identifies $Y_{q+1}$ with $R_q:V^2=U^q-U$. Indeed,
\[
 \left(\frac{U-\theta}{U-\theta^q}\right)^{q+1}-1
 =\frac{(\theta^q-\theta)(U^q-U)}
        {(U-\theta^q)^{q+1}}.
\]
The branch set of $R_q$ is $\mathbb P^1(\mathbb F_q)$, whose
stabilizer is exactly $\operatorname{PGL}_2(\mathbb F_q)$: images
of three distinct rational points determine a unique projective
transformation defined over $\mathbb F_q$.

For $A=\left(\begin{smallmatrix}a&b\\c&d\end{smallmatrix}\right)$
in $\operatorname{GL}_2(\mathbb F_q)$ and $\eta^2=\det A$, the maps
\begin{equation}\label{eq:binomial-Rq-action}
 (U,V)\longmapsto
 \left(\frac{aU+b}{cU+d},
       \frac{\eta V}{(cU+d)^{(q+1)/2}}\right)
\end{equation}
are the two lifts of the corresponding projective transformation.
This follows from
\[
 \left(\frac{aU+b}{cU+d}\right)^q-\frac{aU+b}{cU+d}
 =\frac{\det A\,(U^q-U)}{(cU+d)^{q+1}}.
\]
Their composition agrees with multiplication of pairs in
\eqref{eq:binomial-Rq-group}. Such a pair acts trivially precisely
when $A=\lambda I$ and $\eta=\lambda^{(q+1)/2}$, giving exactly
the denominator of that expression. All automorphisms descend
through the unique hyperelliptic map, so these lifts give the full
group. Its order is
$2|\operatorname{PGL}_2(\mathbb F_q)|=2q(q^2-1)$, as claimed.
\end{proof}

\begin{corollary}\label{cor:quadratic-Cheb-full}
Let $n\geq5$ and assume that the characteristic is zero or prime to
$2n$. The full automorphism group of the smooth curve
$X_{n,2}:y^n=\varphi_2(x)=x^2-2$ is
\[
 \operatorname{Aut}(X_{n,2})\cong
 \begin{cases}
  C_{2n},&n\text{ odd},\\
  \mathcal G_n,&n\text{ even and }n\ne p^e+1
                             \text{ for every }e\geq1,\\
  \mathcal R_q,&n=q+1,\quad q=p^e.
 \end{cases}
\]
These are all the cases of this family with genus at least two.
\end{corollary}

\begin{proof}
Choose $\alpha,\beta\in k$ with $\alpha^n=2$ and $\beta^2=2$.
The substitution $y=\alpha u$, $x=\beta v$ gives
$v^2=u^n+1$, which is isomorphic over $k$ to $Y_n$.
The theorem therefore applies. Its genus is
$\lfloor(n-1)/2\rfloor$, which is at least two exactly when $n\geq5$.
\end{proof}

\section{Artin--Schreier Chebyshev curves}
\label{sec:AS-Chebyshev}

Let $k$ be an algebraically closed field of characteristic $p>2$.  We use
the normalization
\[
 \varphi_0(X)=2,\qquad \varphi_1(X)=X,\qquad
 \varphi_{m+1}(X)=X\varphi_m(X)-\varphi_{m-1}(X).
\]
Thus $\varphi_m(U+U^{-1})=U^m+U^{-m}$.  For $p\nmid m$, let
$A_m$ denote the smooth projective curve with affine equation
\begin{equation}\label{eq:AS-Cheb}
 y^p-y=\varphi_m(x).
\end{equation}
Its genus is $(p-1)(m-1)/2$ \cite[Proposition~3.7.8]{Stichtenoth}.  Write $C_n$ for a cyclic group of order $n$
and $D_n$ for a dihedral group of order $2n$.

\begin{theorem}\label{thm:AS-full-classification}
Suppose $m\geq 3$ and $p\nmid m$.
If $p\geq5$, then
\[
 \operatorname{Aut}_k(A_m)\cong
 \begin{cases}
 C_p\times C_2,&m\text{ even},\\
 D_p,&m\text{ odd}.
 \end{cases}
\]
In characteristic $3$, the same answer holds except for the following five
values of $m$:
\[
\begin{array}{c|c|c}
 m&\operatorname{Aut}_k(A_m)&|\operatorname{Aut}_k(A_m)|\\ \hline
 4&E_{27}\rtimes C_2&54\\
 5&C_3\rtimes C_{10}\cong C_5\times S_3&30\\
 7&D_9&18\\
 8&C_{12}&12\\
 14&C_3\rtimes C_4&12
\end{array}
\]
Here $E_{27}$ is the extraspecial group of order $27$ and exponent $3$.
The indicated involution acts as $-1$ on $E_{27}/Z(E_{27})$ and trivially
on $Z(E_{27})$.  In the last row a generator of $C_4$ acts on $C_3$ by
inversion, so this group has presentation
\[
 \langle u,v\mid u^3=v^4=1,\ vuv^{-1}=u^{-1}\rangle.
\]
Every automorphism in the theorem fixes the unique point above $x=\infty$
and normalizes the Artin--Schreier group
$\langle\rho\rangle$, where $\rho(x,y)=(x,y+1)$.
\end{theorem}

We first show that every automorphism fixes the point above infinity.
We then compute the possible translations and scalings.

\subsection{Reduction to affine substitutions}

For $f\in k[X]$, let $\operatorname{red}(f)$ be its representative modulo
$\{h^p-h:h\in k[X]\}$ with zero constant term and no nonzero term whose
degree is divisible by $p$.  The constant term can be removed because
$k$ is algebraically closed.  Explicitly, if $f=\sum_i a_iX^i$, then
\begin{equation}\label{eq:AS-reduction}
 \operatorname{red}(f)=
 \sum_{\substack{j\geq1\\p\nmid j}}
 \left(\sum_{r\geq0}a_{jp^r}^{p^{-r}}\right)X^j.
\end{equation}
The inverse Frobenius powers are unambiguous in $k$.  Put
$F_m=\operatorname{red}(\varphi_m)$.  Since the coefficients of
$\varphi_m$ lie in $\mathbb F_p$, we have $F_m\in\mathbb F_p[X]$.
Replacing $y$ by $y-h(x)$, and if necessary by a constant translate,
gives an isomorphic model
\begin{equation}\label{eq:AS-reduced-model}
 z^p-z=F_m(x).
\end{equation}

\begin{lemma}\label{lem:AS-full-affine}
If $m\geq3$ and $p\nmid m$, every automorphism of $A_m$ has the form
\begin{equation}\label{eq:AS-affine-form}
 (x,y)\longmapsto (ax+b,cy+h(x)),
 \qquad a\in k^*,\quad b\in k,\quad c\in\mathbb F_p^*,
 \quad h\in k[X],
\end{equation}
where
\begin{equation}\label{eq:AS-affine-condition}
 \varphi_m(ax+b)-c\varphi_m(x)=h(x)^p-h(x).
\end{equation}
Necessarily $c=a^m$ and $\deg h\leq\lfloor(m-1)/p\rfloor$.
Conversely every substitution satisfying
\eqref{eq:AS-affine-condition} defines an automorphism.
\end{lemma}

\begin{proof}
There is one point $P_\infty$ above infinity, and the pole orders of $x$
and $y$ there are $p$ and $m$.  The numerical semigroup generated by
$p$ and $m$ has exactly $(p-1)(m-1)/2$ gaps, so it is the Weierstrass
semigroup at $P_\infty$.

We apply Stichtenoth's theorem \cite[Satz~6]{Stichtenoth1973II}.
We use the formulation in \cite[Theorem~3.1, Proposition~3.3 and
Corollary~3.4]{LehrMatignon2005}, with the arXiv numbering.
It says that every automorphism fixes $P_\infty$ and has the form
\eqref{eq:AS-affine-form}, except for two normal forms.
Up to an affine change of $x$ and an Artin--Schreier change of $y$,
these are
\[
 y^p-y=x^m,\quad m<p,\ m\mid p+1;
 \qquad
 y^p-y=x^{p+1}.
\]
Neither normal form can occur here.  Indeed, when $m<p$ the
Artin--Schreier change is constant.  The coefficient of $x^{m-1}$ in an
affine transform of $x^m$ forces its translation to be zero.  Its
coefficient of $x^{m-2}$ then vanishes, whereas the corresponding
coefficient of $\varphi_m$ is $-m\ne0$.  When $m=p+1$, an affine
transform of $x^{p+1}$ has nonzero terms only in degrees
$p+1,p,1,0$, and an allowed Artin--Schreier change has degree at most
one.  The coefficient of $x^{p-1}$ in $\varphi_{p+1}$ is $-1$, so this
case is excluded as well.

The cited result therefore applies to the full automorphism group.
When $m>p$, the same form follows from Riemann--Roch spaces.
The space $\mathcal L(pP_\infty)$ is spanned by $1,x$.
The space $\mathcal L(mP_\infty)$ is spanned by $y$ and the powers
$x^j$ with $pj\leq m$. When $m<p$ the roles of $x$ and $y$ in the
semigroup are reversed, and one uses instead
\cite[Corollary~3.4]{LehrMatignon2005}, which rests on the normality of
$\langle\rho\rangle$ in the stabilizer of $P_\infty$.  Substitution
in the equation gives \eqref{eq:AS-affine-condition}.  Its leading
coefficient gives $c=a^m$, because a nonconstant polynomial of the form
$h^p-h$ has degree divisible by $p$.  The remaining degree is at most
$m-1$, giving the stated bound on $h$.  The converse follows directly
by substitution, since the displayed transformation is invertible.
\end{proof}

For the remainder of the proof, let
\begin{align}
 S_m&=\{b\in k:\operatorname{red}(F_m(X+b)-F_m(X))=0\},
       \label{eq:AS-S}\\
 J_m&=\{j\geq1:[X^j]F_m\ne0\},\qquad
 e_m=\gcd\bigl(m(p-1),\{m-j:j\in J_m\}\bigr).
       \label{eq:AS-e}
\end{align}
In particular, if $F_m=X^m$, the second definition gives
$e_m=m(p-1)$.  The notation $[X^j]f$ denotes a coefficient.

\Needspace{15\baselineskip}
\begin{lemma}\label{lem:AS-affine-splitting}
The image of $\operatorname{Aut}_k(A_m)$ in the affine group of the
$x$-line is $S_m\rtimes\mu_{e_m}$.  There is an exact sequence
\[
 1\longrightarrow C_p\longrightarrow P_m\longrightarrow S_m
 \longrightarrow1
\]
with central kernel, and
\[
 \operatorname{Aut}_k(A_m)=P_m\rtimes\mu_{e_m},
 \qquad |P_m|=p|S_m|.
\]
Here $S_m$ is a finite $\mathbb F_p$-vector space.  On the reduced
model, the subgroup $\mu_{e_m}$ has the explicit lifts
\begin{equation}\label{eq:AS-scalar-lift}
 (x,z)\longmapsto(ax,a^m z),\qquad a^{e_m}=1.
\end{equation}
The action of such a lift on the kernel $C_p$ is multiplication by
$a^m$.
\end{lemma}

\begin{proof}
The kernel of the action on $k(x)$ is exactly the group
$\langle\rho\rangle\cong C_p$.  A lift of $x\mapsto x+b$ has $c=1$
in \eqref{eq:AS-affine-form}, and therefore commutes with $\rho$.
Such translations form an additive subgroup $S_m$, hence an
$\mathbb F_p$-vector space.  It is finite because the curve has genus
at least two.

The parity identity
$\varphi_m(-X)=(-1)^m\varphi_m(X)$ gives the automorphism
$\iota(x,y)=(-x,(-1)^m y)$.  If $M(X)=aX+b$ belongs to the affine
image, the commutator of $M$ with $X\mapsto-X$ gives the translation
$X\mapsto X-2b$.  Thus $-2b\in S_m$, and hence $b\in S_m$.  Composing
with this translation shows that $X\mapsto aX$ also belongs to the
image.

Since $F_m$ has no positive degree divisible by $p$, the polynomial
$F_m(aX)-a^mF_m(X)$ is reduced.  It is Artin--Schreier trivial if and
only if it is zero.  Together with $a^m\in\mathbb F_p^*$, this says
precisely that $a^{e_m}=1$, and gives the lifts
\eqref{eq:AS-scalar-lift}.  These lifts form a cyclic subgroup and
split the affine multiplier map.  Its kernel is $P_m$, with the exact
sequence asserted in the lemma.
\end{proof}

\subsection{Coefficient calculations}

The coefficients of the Chebyshev polynomial are given by
\begin{equation}\label{eq:AS-Cheb-coeff}
 \begin{aligned}
 \varphi_m(X)=\sum_{j=0}^{\lfloor m/2\rfloor}c_jX^{m-2j},
 &\qquad c_0=1,\\
 c_j=(-1)^j\frac{m}{m-j}\binom{m-j}{j}
 &\qquad(j\geq1).
 \end{aligned}
\end{equation}
The coefficients are integers and are reduced modulo $p$ only after
this expression is evaluated.  In particular,
\[
 c_1=-m,\qquad c_2=\frac{m(m-3)}2,\qquad
 c_3=-\frac{m(m-4)(m-5)}6.
\]
\begin{align*}
 c_4&=\frac{m(m-5)(m-6)(m-7)}{24},\\
 c_5&=-\frac{m(m-6)(m-7)(m-8)(m-9)}{120}.
\end{align*}
For odd $s$, the coefficient of $X^{m-s}$ in
$\varphi_m(X+b)-\varphi_m(X)$ is
\begin{equation}\label{eq:AS-translation-coeff}
 d_s(b)=\sum_{j=0}^{(s-1)/2}
 c_j\binom{m-2j}{s-2j}b^{s-2j},
\end{equation}
provided $m-s\geq0$.

We will use the following coefficient test.  If
$D=h^p-h$ and $\deg D\leq m-1$, then every coefficient $[X^j]D$
with $p\nmid j$ and $j>(m-1)/p$ is zero.  Indeed,
$\deg h\leq(m-1)/p$, so $-h$ contributes nothing in that degree,
and $h^p$ has only degrees divisible by $p$.  Also, a coefficient of
$\varphi_m$ in degree $j>m/p$, with $p\nmid j$, is unchanged by
Artin--Schreier reduction.

\begin{lemma}\label{lem:AS-translations}
For $p\geq5$ one has $S_m=\{0\}$.  In characteristic $3$ the only
nonzero translation spaces are
\[
 S_4=\{b\in k:b^9+b^3+b=0\},\qquad S_7=\mathbb F_3.
\]
Consequently $|S_4|=9$ and $|S_7|=3$.
\end{lemma}

\begin{proof}
Let $b\in S_m$.  We may work with the unreduced polynomial
$\varphi_m$, since an Artin--Schreier change does not alter $S_m$.
If $m\not\equiv1\pmod p$, the coefficient of $X^{m-1}$ in the
translation difference is $mb$.  This degree is prime to $p$ and
greater than $(m-1)/p$, so $b=0$.

Suppose $m\equiv1\pmod p$.  If $p\geq5$, then $m\geq p+1$, and
the coefficient of $X^{m-3}$ is
\[
 \binom m3b^3-m(m-2)b=b.
\]
The degree $m-3$ is prime to $p$ and greater than $(m-1)/p$, giving
$b=0$ again.

Now let $p=3$.  First suppose $m\equiv1\pmod3$ and $m\geq7$.
The same coefficient, now evaluated as an integer before reduction,
is
\begin{equation}\label{eq:AS-first-translation-test}
 d_3(b)=\frac{m-1}{3}b^3+b.
\end{equation}
It must vanish.  In particular, $m\equiv1\pmod9$ gives $b=0$.
When $m\equiv7\pmod9$ it gives $b^3=b$, and when
$m\equiv4\pmod9$ it gives $b^3=-b$.

If $m\equiv7\pmod9$ and $m>7$, formula
\eqref{eq:AS-translation-coeff} gives
\[
 d_5(b)=-b^3=-b.
\]
The degree $m-5$ is prime to $3$ and exceeds $(m-1)/3$; hence $b=0$.
The congruences used here are
$\binom m5=0$, $\binom{m-2}3=1$, and $c_2(m-4)=0$ in
$\mathbb F_3$.

Consider the remaining class $m=4+9r$, and suppose $b\ne0$.
Then $b^3=-b$.  Substitution in
\eqref{eq:AS-translation-coeff} gives
\begin{equation}\label{eq:AS-translation-last-tests}
 d_9(b)=(r+1)b,\qquad d_{11}(b)=(1-r)b
 \quad\text{in characteristic }3.
\end{equation}
Here the integers $r$ are read in $\mathbb F_3$.  These identities
follow directly from Lucas' congruence applied to
\eqref{eq:AS-translation-coeff}: in this residue class
$c_1=-1$, $c_2=-1$, $c_3=0$, $c_4=-1$, and $c_5(m-10)=0$; also
$b^9=b$ and $b^{11}=-b$.  More explicitly, for $d_9$ the only
nonzero contributions are
$\binom m9b^9=rb$ and $c_4(m-8)b=b$.
For $d_{11}$ they are
$c_1\binom{m-2}9b^9=-rb$ and
$c_4\binom{m-8}3b^3=b$.

For $m>13$, the degree $m-9$ is prime to $3$ and exceeds
$(m-1)/3$.  The first identity in
\eqref{eq:AS-translation-last-tests} therefore forces
$r\equiv2\pmod3$.  In this case $m\geq22$, and the degree $m-11$
also satisfies the same degree test.  The second identity then gives
$-b=0$, a contradiction.

The boundary value $m=13$ requires Artin--Schreier reduction in degree
$4$.  Under $b^3=-b$, the coefficients of $X^4$ and $X^{12}$ in
$\varphi_{13}(X+b)-\varphi_{13}(X)$ are respectively $2b$ and $b$.
There is no other degree in the same $3$-power chain, and
$b^{1/3}=-b$.  The reduced coefficient of $X^4$ is therefore
$2b+b^{1/3}=b\ne0$.  Thus this value gives no nonzero translations.

It remains to check $m=4$ and $m=7$.  Their reduced polynomials are
\[
 F_4(X)=X^4-X^2,\qquad F_7(X)=X^7-X^5+X.
\]
For $m=4$,
\[
 F_4(X+b)-F_4(X)
 =bX^3+(b^3+b)X+b^4-b^2.
\]
Its reduced nonconstant part vanishes exactly when
$b^{1/3}+b^3+b=0$, or equivalently $b^9+b^3+b=0$.
This separable additive polynomial has nine roots.  For $m=7$,
\[
 F_7(X+1)-F_7(X)=(X^2+X)^3-(X^2+X)+1.
\]
The constant is Artin--Schreier trivial in $k$, so $1\in S_7$.
Equation \eqref{eq:AS-first-translation-test} already gives
$S_7\subseteq\mathbb F_3$, proving the assertion.
\end{proof}

\begin{lemma}\label{lem:AS-scalars}
If $p\geq5$, then $e_m=2$ for every $m\geq3$ prime to $p$.
If $p=3$, then
\[
 e_m=\begin{cases}
 10,&m=5,\\
 4,&m=8\text{ or }14,\\
 2,&\text{otherwise}.
 \end{cases}
\]
\end{lemma}

\begin{proof}
Let $a\in\mu_{e_m}$.  The parity of the Chebyshev polynomial is
preserved by Artin--Schreier reduction, because $p$ is odd.
Consequently $-1\in\mu_{e_m}$, and it suffices to give upper bounds.

If $m\not\equiv2\pmod p$, then the term $-mX^{m-2}$ is prime to
$p$ in degree and is unchanged by reduction.  Comparing it with the
leading term gives $a^2=1$.

Suppose $p\geq5$ and $m\equiv2\pmod p$.  Since $m\geq3$, we
have $m\geq p+2$.  The coefficient $c_2$ in degree $m-4$ is $-1$
modulo $p$, and this degree is prime to $p$ and exceeds $m/p$.
It follows that $a^4=1$.  The coefficient $c_3$ in degree $m-6$ is
$-2$ modulo $p$; its degree is prime to $p$ and exceeds $m/p$ unless
$(p,m)=(5,7)$.  Outside that boundary case, we obtain $a^6=1$ and
hence $a^2=1$.  In the boundary case the reduced polynomial is
$F_7=X^7-X^3+X$ in characteristic $5$, and its two lower terms give
the same conditions $a^4=a^6=1$.

It remains to consider $p=3$ and $m\equiv2\pmod3$.
For $m\geq8$, the coefficient $c_2=-1$ in degree $m-4$ is
unchanged by reduction, so $a^4=1$.  If $m$ is odd, combining this
with $a^{2m}=1$ gives $a^2=1$.  Suppose $m$ is even and $m>9$.
Writing $t=(m-2)/3$, the integral expression for $c_3$ gives
\[
 c_3=1-t\quad\text{in }\mathbb F_3.
\]
The degree $m-6$ is prime to $3$ and exceeds $m/3$, so this term
gives $a^6=1$ unless $m\equiv5\pmod9$.

In the latter residue class, $c_5=-1$ modulo $3$.  One may verify
this directly by canceling the factor $3$ in $m-8$ against the
denominator $120$ in \eqref{eq:AS-Cheb-coeff}; the remaining
denominator is prime to $3$.  If $m>15$, its degree $m-10$ is prime
to $3$ and exceeds $m/3$, giving $a^{10}=1$.  Combined with $a^4=1$
this again forces $a^2=1$.  Since $m$ is even, the only remaining
value in this residue class is $m=14$.

The three small values are obtained from their reduced polynomials:
\[
\begin{array}{c|c|c}
 m&F_m(X)&e_m\\ \hline
 5&X^5&10\\
 8&X^8-X^4&4\\
 14&X^{14}-X^{10}+X^2&4.
\end{array}
\]
The values of $e_m$ follow either by direct substitution or from
\eqref{eq:AS-e}.
\end{proof}

\subsection{Identification of the groups}

\begin{proof}[Proof of Theorem~\ref{thm:AS-full-classification}]
If $S_m=0$ and $e_m=2$, Lemma~\ref{lem:AS-affine-splitting}
gives $C_p\rtimes C_2$.  The lift of $-1$ is
$\iota(x,z)=(-x,(-1)^m z)$; it centralizes $C_p$ when $m$ is even
and inverts it when $m$ is odd.  This proves every nonexceptional
case, including all cases in characteristic at least five.

In characteristic $3$, the values $m=5,8,14$ have no nontrivial
translations.  Formula \eqref{eq:AS-scalar-lift} gives the action of
the cyclic complement on $C_3$.  For $m=5$ the group is
$C_3\rtimes C_{10}$, with the complement acting through its quotient
$C_2$ by inversion; therefore it is $C_5\times S_3$.  For $m=8$,
the action of $\mu_4$ is trivial, giving $C_3\times C_4\cong C_{12}$.
For $m=14$, a generator $a\in\mu_4$ has $a^{14}=-1$, giving
$C_3\rtimes C_4$ with the action stated in the theorem.

For $m=7$, use the reduced model $z^3-z=X^7-X^5+X$.  Choose
$\beta\in k$ with $\beta^3-\beta=1$.  The automorphism
\[
 \tau(x,z)=(x+1,z+x^2+x+\beta)
\]
has
\[
 \tau^3(x,z)=(x,z-1)=\rho^{-1}(x,z).
\]
Thus the wild subgroup, of order $3|S_7|=9$, is cyclic of order
$9$.  The involution $\iota(x,z)=(-x,-z)$ satisfies
$\iota\tau\iota=\tau^{-1}$, giving $D_9$.

Finally, let $m=4$.  For $b\in S_4$ choose $\beta\in k$ satisfying
$\beta^3-\beta=b^4-b^2$.  The translations lift as
\begin{equation}\label{eq:AS-Heisenberg-lift}
 \tau_{b,\beta}(x,z)
 =(x+b,z+b^{1/3}x+\beta).
\end{equation}
Indeed $b^{1/3}=-b^3-b$ on $S_4$, so substitution gives exactly
the translation difference of $F_4$.  The $27$ transformations
\eqref{eq:AS-Heisenberg-lift} form the wild subgroup $P_4$.
Every element has cube equal to the identity.  Its commutator pairing,
up to the convention for the order of a commutator, is
\[
 \omega(b,c)=b^{1/3}c-c^{1/3}b
            =bc^3-cb^3\in\mathbb F_3.
\]
The second equality follows from the defining equation of $S_4$,
and it also gives $\omega^3=\omega$.  For $b\ne0$, the equation
$\omega(b,c)=0$ is equivalent to $c/b\in\mathbb F_3$.
Consequently this alternating form on the two-dimensional
$\mathbb F_3$-space $S_4$ is nondegenerate.  Hence $P_4$ is the
extraspecial group of order $27$ and exponent $3$.  The involution
$\iota(x,z)=(-x,z)$ induces $b\mapsto-b$ on $S_4$ and fixes the
central deck group, proving the remaining assertion.
\end{proof}

\subsection{The degree-two boundary case}

The hypothesis $m\geq3$ in
Theorem~\ref{thm:AS-full-classification} is essential. The remaining
case of genus at least two is the curve studied by
Roquette~\cite{Roquette1970}, and its automorphisms can move the point
above infinity.

\begin{proposition}\label{prop:AS-degree-two}
Let $p\geq5$. The curve $A_2$ is isomorphic over $k$ to
$R_p:V^2=U^p-U$. Its full group is $\mathcal R_p$, defined in
\eqref{eq:binomial-Rq-group}, of order $2p(p^2-1)$. It fits into the
central exact sequence
\[
 1\longrightarrow C_2\longrightarrow\operatorname{Aut}_k(A_2)
 \longrightarrow\operatorname{PGL}_2(\mathbb F_p)
 \longrightarrow1.
\]
Its explicit transformations are given by
\eqref{eq:binomial-Rq-action} with $q=p$.
\end{proposition}
\begin{proof}
Choose $c\in k$ with $c^p-c=-2$. Setting $U=y-c$ and $V=x$ in
$y^p-y=x^2-2$ gives $V^2=U^p-U$. The assertion is the case $q=p$
of the group calculation in Theorem~\ref{thm:binomial-hyperelliptic}.
\end{proof}

\begin{remark}\label{rem:AS-low-genus}
For $m=1$, the equation $y^p-y=x$ has rational function field $k(y)$,
so its smooth model has genus zero and automorphism group
$\operatorname{PGL}_2(k)$. For $p=3,m=2$ the genus is one and the full
group contains all translations after an origin is chosen. Thus
Theorem~\ref{thm:AS-full-classification} and
Proposition~\ref{prop:AS-degree-two} cover every genus-at-least-two
Artin--Schreier Chebyshev curve in odd characteristic, after applying
Proposition~\ref{prop:Frobenius-reduction}.
\end{remark}
\section{Exceptional deck groups and remaining cases}
\label{sec:remaining}

The following examples show that the deck group need not be normal.
In these cases, the normalizer is a proper subgroup of the full group.

\begin{example}\label{ex:non-normal-quadratic-Cheb}
Let $m=2$ and let $n\geq6$ be even. Since
$\varphi_2(x)=x^2-2$, the change $x=\sqrt2\,v$ and
$y=2^{1/n}u$ gives
\[
 v^2=u^n+1.
\]
The original deck generator becomes $R(u,v)=(\zeta_nu,v)$.
The hyperelliptic involution $\iota(u,v)=(u,-v)$ and the automorphism
\[
 S(u,v)=\left(u^{-1},\frac{v}{u^{n/2}}\right)
\]
satisfy $SRS=\iota R^{-1}$. Since $\iota\notin\langle R\rangle$,
the original deck group $\langle R\rangle$ is not normal. This example
works in characteristic zero and in all characteristics prime to $2n$.
\end{example}

\begin{example}\label{ex:non-normal-Fermat-Cheb}
If $n=m=d\geq4$ and $2d=p^r+1$, the transformation used in
Theorem~\ref{thm:explicit-normalizers} identifies $X_{d,d}$ with a Fermat
curve. Its full automorphism group is
$(C_d\times C_d)\rtimes S_3$, of order $6d^2$, whereas the normalizer of
the distinguished deck group has order $2d^2$. Coordinate permutations
move one coordinate scaling group to another. Thus the normalizer has
index three. The full-group assertion here is the Fermat case of the
plane-curve classification in \cite{TTplane}; for Fermat curves see also
\cite{Tzermias1995} in characteristic zero and \cite{Leopoldt1996} in
positive characteristic.
For example, $(n,m,p)=(6,6,11)$ gives full group order $216$ and
normalizer order $72$.
\end{example}

\begin{example}\label{ex:non-normal-augmented-Fermat}
If $n=m+1$ and $2m+1=p^r$, the transformation used in the augmented
dihedral case gives $v^n=1-z^n$, again a Fermat curve. For $n\geq4$ this
is not a Hermitian degree: $n-1=m=(p^r-1)/2$ is greater than one and not
divisible by $p$. The same Fermat classification
\cite{Leopoldt1996,TTplane} gives full group
order $6n^2$, while the normalizer
of the distinguished deck group has order $2n^2$. Thus allowing infinity
to have the same Kummer exponent as the finite roots gives another
infinite family in which the deck group is not normal.
\end{example}

\subsection*{The remaining Kummer cases}
The table records the cases determined here. It concerns curves of
genus at least two under the characteristic assumptions of Section~4.
\begin{center}
\renewcommand{\arraystretch}{1.2}
\begin{tabular}{p{0.37\textwidth}p{0.55\textwidth}}
\hline
Parameters & Conclusion \\
\hline
$n=2$ & Complete classification in Theorem~\ref{thm:hyp-classification}.\\
$m=2$ & Complete classification in Corollary~\ref{cor:quadratic-Cheb-full}.\\
$m\geq6$, $m>2n$ & The normalizer in Theorem~\ref{thm:explicit-normalizers}
is the full group.\\
$m\geq3$, $n>2m$ & The full group is $\mathcal A_{n,m}$
(Theorem~\ref{thm:full-Cheb-high-n-sharpened}).\\
$n=m$ & The plane-curve classification is taken from \cite{TTplane}.\\
Other parameters & Theorem~\ref{thm:kummer-normalizer} describes the
normalizer.\\
\hline
\end{tabular}
\end{center}
The case $m=5>2n$ necessarily has $n=2$ and is covered by the first row.
The exceptional Fermat models discussed above also have known full groups.
We do not determine the full group in the remaining part of the range
$m\leq2n$, $n\leq2m$. This leads to the following question.

\begin{question}\label{q:middle-range}
Let $3\leq m\leq2n$ and $3\leq n\leq2m$ with $n\ne m$. Is the
distinguished deck group $N$ normal in $\operatorname{Aut}(X_{n,m})$
outside the Fermat family of Example~\ref{ex:non-normal-augmented-Fermat}?
Equivalently, does $\operatorname{Aut}(X_{n,m})$ then coincide with the
normalizer $M$ described in Theorem~\ref{thm:kummer-normalizer}?
\end{question}

\section*{Acknowledgements}

The second author acknowledges partial financial support from CNPq,
grant no.~302774/2025-4, and FAPESP, grant no.~2024/00923-6.

\end{document}